\documentclass[11pt]{amsart}
\usepackage{amsmath,amssymb,amsthm}
\usepackage[margin=1.15in]{geometry}
\usepackage{graphicx}
\usepackage[colorlinks,linkcolor=blue,citecolor=red]{hyperref}

\newtheorem{thm}{Theorem}[section]
\newtheorem{lema}[thm]{Lemma}
\newtheorem{prop}[thm]{Proposition}
\newtheorem{cor}[thm]{Corollary}
\theoremstyle{definition}

\newtheorem{exam}[thm]{Example}
\theoremstyle{remark}
\newtheorem{rem}[thm]{Remark}
\numberwithin{equation}{section}

\newcommand{\R}{\mathbb R}

\newcommand{\Sn}{\mathbb S^{n-1}}
\newcommand{\Hn}{\mathcal H^{n-1}}
\newcommand{\ve}{\varepsilon}
\newcommand{\A}{\mathcal A}
\newcommand{\E}{\mathcal E}
\newcommand{\intr}{\int_{\R^n}}

\begin{document}

\title[Localization and homogenization of fractional energies]{Localization and homogenization of anisotropic fractional energies}

\author[J. Fern\'andez Bonder]{Juli\'an Fern\'andez Bonder}
\address{Departamento de Matem\'atica, Facultad de Ciencias Exactas y Naturales, Universidad de Buenos Aires, Buenos Aires, Argentina, and Instituto de C\'alculo UBA-CONICET, Buenos Aires, Argentina}
\email{jfbonder@dm.uba.ar}

\author[A.   Salort]{Ariel  Salort}
\address{Departamento de Matem\'aticas Aplicadas y Ciencia de Datos, Universidad San Pablo-CEU, CEU Universities, Juli\'an Romea 2, 28003 Madrid, Spain}
\email{ariel.salort@ceu.es}


\begin{abstract}
We study the joint limit of localization ($s\uparrow1$) and periodic homogenization ($\ve\downarrow0$) for anisotropic fractional energies with $p$-growth and oscillating kernels. We prove that the $\Gamma$-limit is governed by the single parameter $\theta=\lim\ve^{p(1-s)}$, and that its density is the convex combination $\theta\A^*+(1-\theta)\bar\A$ of the homogenized and the averaged densities. This shows that the picture recently obtained for quadratic energies also holds for every $p$ and for anisotropic kernels, with no regularity in the oscillating variable; the approach is different, and is based on convex duality. A nonlocal version of the duality argument gives the limit as $s\uparrow1$ for a fixed kernel under weak assumptions, correcting a previous proof by the authors.
\end{abstract}

\maketitle

\section{Introduction}

A celebrated result of Bourgain, Brezis and Mironescu \cite{BBM} states that, for $u\in L^p(\R^n)$,
$$
\lim_{s\uparrow1}(1-s)\iint_{\R^n\times\R^n}\frac{|u(x)-u(y)|^p}{|x-y|^{n+sp}}\,dx\,dy=\frac{K_{n,p}}p\intr|\nabla u|^p\,dx,\qquad K_{n,p}=\int_{\Sn}|\omega\cdot e_1|^p\,d\Hn ,
$$
so that the nonlocal energies \emph{localize} as $s\uparrow1$; the convergence also holds in the sense of $\Gamma$-convergence \cite{Ponce}. Homogenization is another mechanism which replaces a complicated energy by an effective local one: if the coefficients oscillate at a small scale $\ve$, the energies converge as $\ve\downarrow0$ to an energy with constant coefficients, which are not the averages of the original ones. In this paper we study what happens when both mechanisms act at the same time.

We consider the anisotropic fractional energies
\begin{equation}\label{intro.E}
\E_{\ve,s}(u)=(1-s)\iint_{\R^n\times\R^n}m\Big(\frac x\ve,h\Big)\frac{|u(x)-u(x-h)|^p}{|h|^{n+sp}}\,dx\,dh ,\qquad1<p<\infty,
\end{equation}
where the kernel $m(y,h)$ is periodic in $y$, bounded and bounded away from zero, and $m(y,r\omega)$ converges as $r\to0$ to a profile $a(y,\omega)$, which may depend on the direction $\omega\in\Sn$. For $\ve=1$ these energies were introduced in \cite{FBS-Israel}, where it was shown that they converge as $s\uparrow1$ to the local energy $\intr\A(x,\nabla u)\,dx$ with
$$
\A(y,\xi)=\frac1p\int_{\Sn}a(y,\omega)|\xi\cdot\omega|^p\,d\Hn(\omega).
$$
In the same paper it was observed that, for the oscillating kernels in \eqref{intro.E}, the two iterated limits are different. If $s\uparrow1$ first and then $\ve\downarrow0$, one obtains the homogenized density
$$
\A^*(\xi)=\inf\Big\{\int_Q\A(y,\xi+\nabla v(y))\,dy\colon v\in W^{1,p}_{\rm per}(Q)\Big\}
$$
given by the cell problem; if $\ve\downarrow0$ first (see \cite{FBRS}) and then $s\uparrow1$, one obtains the average $\bar\A(\xi)=\int_Q\A(y,\xi)\,dy$, where $Q=(0,1)^n$ is the unit cube. Hence the limit as $(s,\ve)\to(1,0)$ must depend on the relation between the two parameters.

Writing the kernel as a function of the point $x$ and of the increment $h$, rather than of the pair $(x,y)$, was the point of view adopted in \cite{FBS-Israel}, and it is essential in what follows: the behavior of $m(x,h)$ as $h\to0$ along each direction is what determines the limit as $s\uparrow1$, and making the kernel oscillate in $x$ alone leads, as $s\uparrow1$, to the homogenization problem for an anisotropic medium (for $p=2$, for an arbitrary periodic elliptic matrix).

\begin{figure}[ht]
\centering
\includegraphics[width=\textwidth]{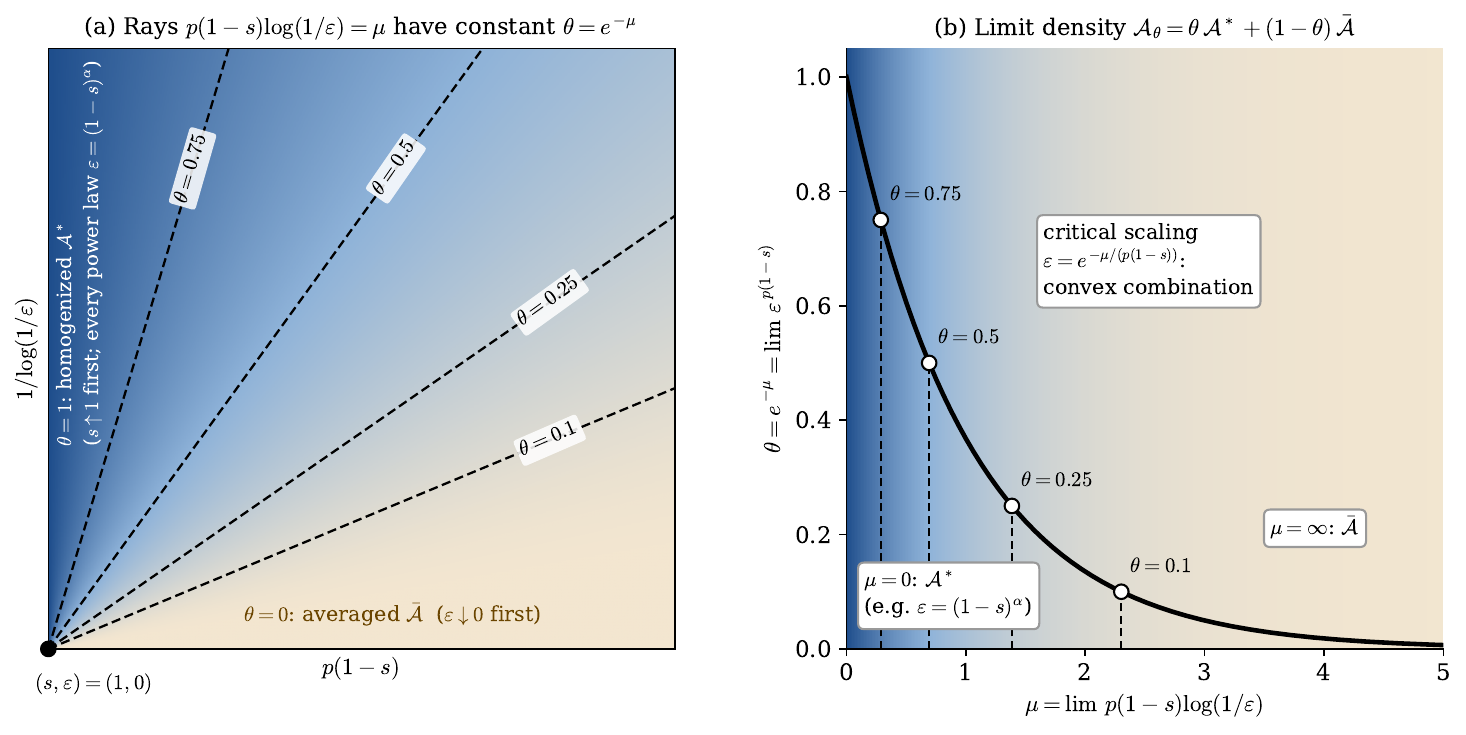}
\caption{Regimes of the joint limit (Theorem \ref{teo.main}). (a) Near the limit point $(s,\ve)=(1,0)$, in the variables $p(1-s)$ and $1/\log(1/\ve)$, the parameter $\theta=\lim\ve^{p(1-s)}$ is constant along the rays $p(1-s)\log(1/\ve)=\mu$, where $\theta=e^{-\mu}$ (that is, $\mu=p\lambda$ with $\lambda=\lim(1-s)\log(1/\ve)$). Paths tangent to the vertical axis, such as all power laws $\ve=(1-s)^\alpha$, give $\theta=1$, and paths tangent to the horizontal axis give $\theta=0$. (b) The weight $\theta=e^{-\mu}$ of $\A^*$ in the limit density $\A_\theta=\theta\A^*+(1-\theta)\bar\A$: the transition from the homogenized to the averaged density takes place at the scale $\ve=e^{-\mu/(p(1-s))}$.}\label{fig.reg}
\end{figure}

\medskip
\noindent\textbf{The joint limit.} Our main result (Theorem \ref{teo.main}) is that the joint limit is governed by the single parameter
$$
\theta=\lim\ve^{p(1-s)}\in[0,1],
$$
and that the energies $\E_{\ve,s}$ $\Gamma$-converge in $L^p(\R^n)$ to the local energy with density
$$
\A_\theta=\theta\,\A^*+(1-\theta)\,\bar\A ,
$$
together with the compactness in $L^p_{\rm loc}(\R^n)$ of bounded sequences with bounded energy and the convergence of minimizers. The reason is that the radial measure $(1-s)\,p\,r^{p(1-s)-1}dr$ gives mass $\rho^{p(1-s)}$ to the scales $r<\rho$: the scales $|h|\ll\ve$, which see the local oscillating problem, carry a fraction $\theta$ of the energy, the scales $|h|\gg\ve$, which see the averaged kernel, carry a fraction $1-\theta$, and the scales $|h|\sim\ve$ are negligible. In particular no nonlocal cell problem appears, the set of possible limits is a segment, and the limit exists along a family $(s,\ve)$ if and only if $\ve^{p(1-s)}$ converges (unless $\A^*=\bar\A$). Note that $\theta=1$ for every power law $\ve=(1-s)^\alpha$; the transition takes place at the exponentially small scale $\ve=e^{-\lambda/(1-s)}$, for which $\theta=e^{-p\lambda}$. See Figure \ref{fig.reg}.

\medskip
\noindent\textbf{Related results.} To our knowledge, the interaction between the Bourgain--Brezis--Mironescu limit and homogenization for fractional energies was first considered in \cite{FBS-Israel}, for the energies \eqref{intro.E} and every $p$: as recalled above, it was observed there that the two iterated limits are different, which is the starting point of the present paper. 

Later, and independently, the problem has been studied for quadratic energies with isotropic kernels. Braides, Brusca and Donati \cite{BBD}   considered energies $(1-s)\iint_{\Omega\times\Omega}k(\frac x\ve)\frac{|u(x)-u(y)|^2}{|x-y|^{n+2s}}dxdy$, with $k$ continuous and periodic, and proved that they $\Gamma$-converge to the homogenized Dirichlet energy when $1-s\ll\ve^2$; they indicate the extension to $p$-Gagliardo seminorms when $1-s\ll\ve^p$. 

Very recently, Braides, Mannella and Schiavoni Piazza \cite{BMS} gave the complete description for the energies $(1-s)\iint_{\Omega\times\Omega}k(\frac x\ve,\frac y\ve)\frac{|u(x)-u(y)|^2}{|x-y|^{n+2s}}dxdy$ on a bounded domain, with $k$ continuous and periodic in both variables: the $\Gamma$-limit is the convex combination of the averaged and the homogenized energies with weights $1-\lambda$ and $\lambda=\lim\ve^{2(1-s)}$. For $p=2$ and an isotropic coefficient depending only on $x$, the limit density is the one of Theorem \ref{teo.main}, with $\lambda=\theta$; neither result contains the other, since in \cite{BMS} the energy is defined on $\Omega\times\Omega$ and the coefficient may depend on both variables, while here the energy is defined on $\R^n\times\R^n$.

The description of the kernel in terms of $(x,h)$ makes a difference: even for $p=2$ the two classes of kernels are different. The kernels $m(\frac x\ve,x-y)$ considered here depend on the direction of $x-y$, so that for fixed $\ve$ the energies localize, as $s\uparrow1$, to $\intr A(\frac x\ve)\nabla u\cdot\nabla u\,dx$, where $A$ is an arbitrary symmetric, uniformly elliptic, periodic matrix: Theorem \ref{teo.main} describes the nonlocal approximation of the homogenization of anisotropic media (Corollary \ref{cor.matrices}). For the coefficients $k(\frac x\ve,\frac y\ve)$ of \cite{BMS} the energies localize to a multiple of $\int_\Omega k(\frac x\ve,\frac x\ve)|\nabla u|^2dx$, with a scalar conductivity. As families in $\ve$, the two classes only share the isotropic coefficients $k(\frac x\ve)$ (up to the different domains of integration). 
In \cite{BMS} it is expected that the same picture holds for $p$-Gagliardo seminorms and for more general functionals. Theorem \ref{teo.main} shows that the same picture does hold for the energies \eqref{intro.E}: for every $1<p<\infty$ and for anisotropic kernels whose profile depends on the direction and is only Borel measurable in the oscillating variable. The situation is different for energies of convolution type, where the kernel is integrable and concentrates at a scale $\delta$: there a general theory is available \cite{AABPT}, the critical scaling is $\ve\sim\delta$, and at the critical scaling the limit is given by a nonlocal cell problem \cite{AABPT} (see \cite{Brusca} for the analysis of the different regimes), instead of the coexistence of the two extreme regimes.

\medskip
\noindent\textbf{The approach.} Our proof is different from the ones in \cite{BBD,BMS}. There the short-range interactions are treated by a discretization and averaging of the functions, introduced in \cite{BBD}, in which the uniform continuity of the coefficient is used to factor it out of the integrals (in \cite[Remark 3.1]{BMS} it is stated that the continuity is not necessary for coefficients depending only on $x$), and in \cite{BMS} the long-range interactions are treated by a blow-up argument and a periodization. We do not discretize and we do not use blow-up. Instead, we rely on two elementary estimates for the energy in a band of scales and on convex duality.

For the large scales the key observation is that, in the band $\{|h|>\rho\}$, replacing $u$ by its mollification at any scale $\delta\le\rho$ changes the energy by $O(1-s)$ (Lemma \ref{lema.moll}); this allows us to average the kernel. For the small scales we mollify at a scale $\delta\ll\ve$ using H\"older's inequality instead of Jensen's (Lemma \ref{lema.holder}), which replaces the weight by a local harmonic-type mean and reduces the problem to a local one with slightly perturbed coefficients; we then use a proof by duality of the lower bound in periodic homogenization, which is stable under such perturbations. Since convex duality is not restricted to quadratic energies, the argument works for every $p$, with the nonlinear cell problem defining $\A^*$, and since the coefficients are never evaluated along the functions --- the weight is first replaced by a local mean, and the homogenization step is done by duality --- no regularity in the oscillating variable is needed and the direction-dependent profile causes no additional difficulty. The lower bounds hold for sequences which are bounded in $L^p(\R^n)$ and converge only in $L^p_{\rm loc}(\R^n)$, and the errors in the reduction steps are explicit (Section \ref{sec.com}). For the joint limit we assume that $m(y,r\omega)$ converges to $a(y,\omega)$ with a rate, uniformly in $y$ (condition \eqref{H3}); we do not know whether the integrated condition \eqref{K} below suffices. Coefficients of the form $k(\frac x\ve,\frac y\ve)$, which oscillate in both variables, are not of the form \eqref{intro.E} and are not covered by our result (see Section \ref{sec.com}).

\medskip
\noindent\textbf{A fixed kernel.} The same duality argument, applied directly at the nonlocal level, gives the limit as $s\uparrow1$ for a fixed kernel under weak assumptions. In Theorem \ref{teo.fijo} we prove that the $\Gamma$-convergence holds as soon as $m(y,t\omega)$ converges to its profile $a(y,\omega)$ in $L^1_{\rm loc}$ in $y$ and in mean in $\omega$, with no rate and no uniformity in $y$ (condition \eqref{K}). This covers, for instance, the symmetric kernels $\frac12(c(x)+c(y))$, $\sqrt{c(x)c(y)}$ or $c(x)c(y)$ with $c$ measurable and bounded between two positive constants. For weights of the form $f(x)f(y)$ the $\Gamma$-convergence as $s\uparrow1$ was recently proved in \cite{KSS} when $f$ is Lipschitz continuous, for $p\ge1$ and also for weights depending on $s$ which converge uniformly; for $p>1$ and weights bounded away from zero, Theorem \ref{teo.fijo} shows that no regularity of $f$ is needed. The lower bound holds for sequences converging only in $L^1_{\rm loc}$: the test functions are smooth in $x$, so that the coefficient only enters through an error term and the functions are never regularized. This extends the result of \cite{FBS-Israel}, in which the convergence $m(x,r\omega)\to a(x,\omega)$ was required with a rate, and it repairs the proof of the lower bound given there, which relied on a lemma that is not correct as stated (see Remark \ref{rem.israel}). The same proof applies to kernels depending on $s$ (Remark \ref{rem.ms}).

\medskip
\noindent\textbf{Organization.} Section \ref{sec.tools} contains the setting and the tools. Section \ref{sec.fijo} deals with a fixed kernel. In Section \ref{sec.conjunto} we prove the result on the joint limit, and Section \ref{sec.com} contains some further comments.

\section{Setting and tools}\label{sec.tools}

\subsection{Kernels and energies}
Let $n\ge1$, $1<p<\infty$ and $Q=(0,1)^n$. Throughout the paper a \emph{kernel} is a measurable function $m\colon\R^n\times\R^n\to\R$ such that
\begin{equation}\label{H1}\tag{$H_1$}
0<m_-\le m(y,h)\le m_+<\infty ,
\end{equation}
and a \emph{profile} is a measurable function $a\colon\R^n\times\Sn\to[m_-,m_+]$, which describes the behavior of $m(y,h)$ as $h\to0$ in the direction $\omega=\frac h{|h|}$. We use two conditions relating $m$ and $a$. The first one is integrated in $y$ and $\omega$ and has no rate: for every $R>0$,
\begin{equation}\label{K}\tag{$K$}
\varpi_R(r):=\operatorname*{ess\,sup}_{0<t<r}\int_{B_R}\frac1{n\omega_n}\int_{\Sn}\big|m(y,t\omega)-a(y,\omega)\big|\,d\Hn(\omega)\,dy\ \longrightarrow\ 0\qquad\text{as }r\to0^+ .
\end{equation}
(The function $(y,t,\omega)\mapsto m(y,t\omega)$ is measurable, so the integral is defined for a.e.\ $t$.) The second one is pointwise and has a rate, uniform in $y$: there is $C_0>0$ such that
\begin{equation}\label{H3}\tag{$H_3$}
|m(y,r\omega)-a(y,\omega)|\le C_0\,r\qquad\text{for all }y\in\R^n,\ \omega\in\Sn,\ 0<r<1 .
\end{equation}
Clearly \eqref{H3} implies \eqref{K}, with $\varpi_R(r)\le C_0|B_R|\,r$. No symmetry of $m$ and no regularity of $m$ or $a$ in the variable $y$ is assumed. For a fixed kernel (Section \ref{sec.fijo}) we only need \eqref{K}; for the joint limit (Section \ref{sec.conjunto}) we assume \eqref{H3}, and that $m$ is Borel measurable and $Q$-periodic in $y$.

For $\ve>0$, $s\in(0,1)$, $0\le\rho_1<\rho_2\le\infty$ and $u\in L^1_{\rm loc}(\R^n)$ we write
$$
\E_{\ve,s}(u;\rho_1,\rho_2)=(1-s)\int_{\{\rho_1<|h|<\rho_2\}}\intr m\Big(\frac x\ve,h\Big)\frac{|u(x)-u(x-h)|^p}{|h|^{n+sp}}\,dx\,dh\in[0,\infty],
$$
$\E_{\ve,s}(u)=\E_{\ve,s}(u;0,\infty)$, and $\E_s=\E_{1,s}$ for the energy of the kernel $m$ itself. The limit density is
$$
\A(y,\xi)=\frac1p\int_{\Sn}a(y,\omega)|\xi\cdot\omega|^p\,d\Hn(\omega).
$$
With this normalization (without the factor $\frac1p$ in front of the energy used in \cite{FBS-Israel}) the limit of $\E_s(u)$ as $s\uparrow1$ is exactly $\intr\A(x,\nabla u)$, since $(1-s)\int_0^1r^{p(1-s)-1}dr=\frac1p$. Since $|\xi\cdot\omega|^p$ is even in $\omega$, only the even part of $a(y,\cdot)$ matters.

We also use the notation $\tau_hu(x)=u(x+h)$ for increments of $u$, $[u]_{s,p}^p=\iint\frac{|u(x)-u(x-h)|^p}{|h|^{n+sp}}dx\,dh$ to denote the $(s,p)$-Gagliardo seminorm, $\omega_n=|B_1|$ and $K_{n,p}=\int_{\Sn}|\omega\cdot e_1|^pd\Hn$. 

\noindent Let $\varrho\in C^\infty_c(B_1)$ be a standard mollifier, $\varrho_\delta(x)=\delta^{-n}\varrho(x/\delta)$, $u_\delta=u*\varrho_\delta$, $c_\varrho=\|\nabla\varrho\|_{L^1}$ and $c_\varrho'=\|D^2\varrho\|_{L^1}$.

\subsection{Translations and mollification}

\begin{lema}[Translation estimate, uniform in $s$]\label{lema.transl}
For every $u\in L^p(\R^n)$, $s\in(0,1)$ and $h\in\R^n$,
$$
\|\tau_hu-u\|_p^p\le C_{n,p}\,|h|^{sp}\,(1-s)[u]_{s,p}^p,\qquad C_{n,p}=\frac{p\,2^{2p+1}4^n}{\omega_n}.
$$
\end{lema}

\begin{proof}
Let $\Psi(r)=r^{-p}\int_{\Sn}\|\tau_{r\omega}u-u\|_p^p\,d\Hn(\omega)$, so that $[u]_{s,p}^p=\int_0^\infty\Psi(r)r^{p(1-s)}\frac{dr}r$. Since
$$
\frac{|u(x+2r\omega)-u(x)|}{2r}\le\frac12\Big(\frac{|u(x+2r\omega)-u(x+r\omega)|}{r}+\frac{|u(x+r\omega)-u(x)|}{r}\Big),
$$
convexity and translation invariance give $\Psi(2r)\le\Psi(r)$ for every $r>0$. Fix $r_0>0$ and let $\varphi_j=\int_{2^{-j-1}r_0}^{2^{-j}r_0}\Psi(r)\frac{dr}r$. The change of variables $r=2^{-j}r'$ gives $\varphi_j\ge\varphi_0$ for $j\ge0$. Hence
$$
\int_0^{r_0}\Psi(r)r^{p(1-s)}\frac{dr}r\ge\sum_{j\ge0}(2^{-j-1}r_0)^{p(1-s)}\varphi_j\ge\frac{r_0^{p(1-s)}}{2^{p(1-s)}-1}\,\varphi_0 .
$$
Since $2^{p(1-s)}-1\le p(1-s)2^p$ and $\Psi(r)\ge r_0^{-p}\int_{\Sn}\|\tau_{r\omega}u-u\|_p^p\,d\Hn$ for $r\le r_0$, we get
$$
\int_{r_0/2}^{r_0}\int_{\Sn}\|\tau_{r\omega}u-u\|_p^p\,d\Hn\frac{dr}r\le r_0^p\varphi_0\le p2^p\,r_0^{sp}\,(1-s)[u]_{s,p}^p .
$$
Applying this also with $r_0/2$ in place of $r_0$,
\begin{equation}\label{anillo}
\int_{r_0/4}^{r_0}\int_{\Sn}\|\tau_{r\omega}u-u\|_p^p\,d\Hn\frac{dr}r\le p2^{p+1}\,r_0^{sp}\,(1-s)[u]_{s,p}^p .
\end{equation}
Now let $|h|=r_0$ and $B=B(h/2,r_0/4)$. For $z\in B$, $\|\tau_hu-u\|_p^p\le2^{p-1}(\|\tau_{h-z}u-u\|_p^p+\|\tau_zu-u\|_p^p)$. Since $z\mapsto h-z$ maps $B$ onto itself, averaging over $B$ gives $\|\tau_hu-u\|_p^p\le2^p|B|^{-1}\int_B\|\tau_zu-u\|_p^p\,dz$. Since $B\subset\{r_0/4\le|z|\le r_0\}$, in polar coordinates
\begin{align*}
\frac1{|B|}\int_B\|\tau_zu-u\|_p^p\,dz&\le\frac{4^n}{\omega_nr_0^n}\int_{r_0/4}^{r_0}\int_{\Sn}\|\tau_{r\omega}u-u\|_p^p\,r^n\,d\Hn\frac{dr}r\\
&\le\frac{4^n}{\omega_n}\int_{r_0/4}^{r_0}\int_{\Sn}\|\tau_{r\omega}u-u\|_p^p\,d\Hn\frac{dr}r,
\end{align*}
and the result follows from \eqref{anillo}.
\end{proof}

By \eqref{H1}, for every $\ve>0$,
\begin{equation}\label{coerc}
(1-s)[u]_{s,p}^p\le\frac1{m_-}\,\E_{\ve,s}(u).
\end{equation}

\begin{lema}[Mollification]\label{lema.moll}
Let $u\in L^p(\R^n)$, $s\in(0,1)$, $\delta>0$ and $\Lambda=(1-s)[u]_{s,p}^p$. Then
\begin{enumerate}
\item $\|u_\delta-u\|_p^p\le C_{n,p}\,\delta^{sp}\Lambda$;
\item $\|\tau_hu_\delta-u_\delta\|_p\le\|\tau_hu-u\|_p$ and $\|\tau_h\nabla u_\delta-\nabla u_\delta\|_p\le\frac{c_\varrho}\delta\|\tau_hu-u\|_p$ for every $h$;
\item for every $\rho\ge\delta$,
$$
(1-s)\int_{\{|h|>\rho\}}\intr\frac{|(u-u_\delta)(x)-(u-u_\delta)(x-h)|^p}{|h|^{n+sp}}\,dx\,dh\le\frac{2^pn\omega_nC_{n,p}}{sp}\,(1-s)\,\Lambda .
$$
\end{enumerate}
\end{lema}

\begin{proof}
(1) follows from $u_\delta-u=\int\varrho_\delta(z)(\tau_{-z}u-u)\,dz$, Minkowski's inequality and Lemma \ref{lema.transl}. (2) follows from $\tau_hu_\delta-u_\delta=(\tau_hu-u)*\varrho_\delta$, $\tau_h\nabla u_\delta-\nabla u_\delta=(\tau_hu-u)*\nabla\varrho_\delta$ and Young's inequality, since $\|\nabla\varrho_\delta\|_1=c_\varrho/\delta$. For (3), let $w=u-u_\delta$. The left-hand side is at most
$$
(1-s)\,2^p\|w\|_p^p\int_{\{|h|>\rho\}}\frac{dh}{|h|^{n+sp}}=(1-s)\,2^p\|w\|_p^p\,\frac{n\omega_n}{sp}\,\rho^{-sp},
$$
and by (1), $\|w\|_p^p\rho^{-sp}\le C_{n,p}(\delta/\rho)^{sp}\Lambda\le C_{n,p}\Lambda$.
\end{proof}

\begin{rem}
Item (3) is a key point: in the band $\{|h|>\rho\}$, replacing $u$ by its mollification at any scale $\delta\le\rho$ changes the energy by a quantity which is $O(1-s)$ times the total energy.
\end{rem}

\begin{lema}[Uniform upper bound]\label{lema.cota}
For every $v\in W^{1,p}(\R^n)$, $\ve>0$, $s\in(0,1)$, $0\le\rho_1<\rho_2<\infty$ and $\rho>0$,
$$
\E_{\ve,s}(v;\rho_1,\rho_2)\le c_2\,\big(\rho_2^{p(1-s)}-\rho_1^{p(1-s)}\big)\|\nabla v\|_p^p,\qquad\E_{\ve,s}(v;\rho,\infty)\le\frac{2^pn\omega_nm_+}{sp}(1-s)\rho^{-sp}\|v\|_p^p ,
$$
where $c_2=\frac{m_+K_{n,p}}p$. In particular, $\E_{\ve,s}(\cdot\,;\rho_1,\rho_2)^{1/p}$ is a seminorm on $W^{1,p}(\R^n)$ and, for $s\ge\frac12$, it is bounded by $C\|\cdot\|_{W^{1,p}}$, with $C$ independent of $\ve$, $s$, $\rho_1$ and $\rho_2\le\infty$.
\end{lema}

\begin{proof}
Since $\|\tau_{-r\omega}v-v\|_p^p\le r^p\intr|\nabla v\cdot\omega|^pdx$ and $\int_{\Sn}|\xi\cdot\omega|^pd\Hn=K_{n,p}|\xi|^p$,
\begin{align*}
\E_{\ve,s}(v;\rho_1,\rho_2)&\le m_+(1-s)\int_{\rho_1}^{\rho_2}r^{p(1-s)-1}dr\int_{\Sn}\intr|\nabla v\cdot\omega|^pdx\,d\Hn\\
&=c_2\big(\rho_2^{p(1-s)}-\rho_1^{p(1-s)}\big)\|\nabla v\|_p^p .
\end{align*}
The second bound follows from $\|\tau_{-h}v-v\|_p^p\le2^p\|v\|_p^p$. For the last assertion, note that the integrand is nonnegative, so $\E_{\ve,s}(v;\rho_1,\rho_2)\le\E_{\ve,s}(v;0,\infty)=\E_{\ve,s}(v;0,1)+\E_{\ve,s}(v;1,\infty)$. By the first bound with $(\rho_1,\rho_2)=(0,1)$ and the second one with $\rho=1$, for $s\ge\frac12$,
$$
\E_{\ve,s}(v;\rho_1,\rho_2)\le c_2\|\nabla v\|_p^p+\frac{2^pn\omega_nm_+}{p}\|v\|_p^p\le C^p\|v\|_{W^{1,p}}^p,
$$
where $C^p=\max\Big\{c_2,\frac{2^pn\omega_nm_+}p\Big\}$. This concludes the proof.
\end{proof}

\begin{lema}[Derivatives of the mollification]\label{lema.der}
Let $u\in L^p(\R^n)$, $s\in(0,1)$, $\delta>0$ and define $\Lambda=(1-s)[u]_{s,p}^p$. Then
$$
\|\nabla u_\delta\|_p\le c_\varrho\,(C_{n,p}\Lambda)^{1/p}\,\delta^{s-1},\qquad\|D^2u_\delta\|_p\le c_\varrho'\,(C_{n,p}\Lambda)^{1/p}\,\delta^{s-2}.
$$
\end{lema}

\begin{proof}
For $|\alpha|\ge1$, $\int\partial^\alpha\varrho_\delta=0$, so $\partial^\alpha u_\delta(x)=\int\partial^\alpha\varrho_\delta(z)(u(x-z)-u(x))\,dz$, and Minkowski's inequality and Lemma \ref{lema.transl} give $\|\partial^\alpha u_\delta\|_p\le\|\partial^\alpha\varrho_\delta\|_1(C_{n,p}\delta^{sp}\Lambda)^{1/p}$.
\end{proof}

\begin{lema}[Weighted mollification]\label{lema.holder}
Let $b\colon\R^n\to[m_-,m_+]$ be measurable, $\delta>0$, and let
$$
b^\delta(x)=\Big(\int\varrho_\delta(z)\,b(x-z)^{-\frac1{p-1}}\,dz\Big)^{-(p-1)} .
$$
Then $m_-\le b^\delta\le m_+$ and, for every $u\in L^p(\R^n)$ and $h\in\R^n$,
$$
\intr b^\delta(x)\,|u_\delta(x)-u_\delta(x-h)|^p\,dx\le\intr b(x)\,|u(x)-u(x-h)|^p\,dx .
$$
\end{lema}

\begin{proof}
Let $d(x)=u(x)-u(x-h)$, so that $u_\delta(x)-u_\delta(x-h)=\int\varrho_\delta(z)d(x-z)\,dz$. Writing $\varrho_\delta d=(\varrho_\delta^{1/p'}b^{-1/p})(\varrho_\delta^{1/p}b^{1/p}d)$, H\"older's inequality gives
$$
|u_\delta(x)-u_\delta(x-h)|^p\le\Big(\int\varrho_\delta(z)b(x-z)^{-\frac{p'}p}dz\Big)^{p-1}\int\varrho_\delta(z)\,b(x-z)\,|d(x-z)|^p\,dz ,
$$
and $\frac{p'}p=\frac1{p-1}$. Multiplying by $b^\delta(x)$ and integrating in $x$, Fubini's theorem gives the result.
\end{proof}

\begin{rem}
With Jensen's inequality instead of H\"older's one would need to compare $b(x-z)$ with $b(x)$ for $|z|\le\delta$, that is, a modulus of continuity. Lemma \ref{lema.holder} avoids this: the price is that the weight $b$ is replaced by its local harmonic-type mean $b^\delta$, which converges to $b$ a.e.\ as $\delta\to0$.
\end{rem}

\subsection{Densities}
For $b\colon\R^n\times\Sn\to[m_-,m_+]$ measurable we write 
$$
\mathcal B(y,\xi)=\frac1p\int_{\Sn}b(y,\omega)|\xi\cdot\omega|^pd\Hn.
$$ 
Every such function is convex and $C^1$ in $\xi$, and satisfies
\begin{equation}\label{crec}
c_1|\xi|^p\le\mathcal B(y,\xi)\le c_2|\xi|^p,\qquad|\nabla_\xi\mathcal B(y,\xi)|\le m_+K_{n,p}\,|\xi|^{p-1},
\end{equation}
with $c_1=\frac{m_-K_{n,p}}p$, $c_2=\frac{m_+K_{n,p}}p$ (for the last bound, $\nabla_\xi\mathcal B\cdot e=\int b|\xi\cdot\omega|^{p-2}(\xi\cdot\omega)(\omega\cdot e)\,d\Hn$ and H\"older's inequality on $\Sn$). Its convex conjugate 
$$\mathcal B^c(y,\sigma)=\sup_\xi\{\sigma\cdot\xi-\mathcal B(y,\xi)\}$$ 
satisfies $0\le\mathcal B^c(y,\sigma)\le c_3|\sigma|^{p'}$ with $c_3$ depending only on $c_1$ and $p$. These facts apply in particular to $\A$.

\section{The limit for a fixed kernel}\label{sec.fijo}

In this section $m$ is a kernel satisfying \eqref{H1} and \eqref{K} for some profile $a$, and $\E_s=\E_{1,s}$.

\begin{thm}\label{teo.fijo}
Assume \eqref{H1} and \eqref{K}. Then:
\begin{enumerate}
\item[(i)] for every $u\in W^{1,p}(\R^n)$, $\displaystyle\lim_{s\uparrow1}\E_s(u)=\intr\A(x,\nabla u)\,dx$;
\item[(ii)] if $s_k\uparrow1$, $u_k\to u$ in $L^1_{\rm loc}(\R^n)$ and $L:=\liminf_k\E_{s_k}(u_k)<\infty$, then $u\in W^{1,1}_{\rm loc}(\R^n)$, $\nabla u\in L^p(\R^n;\R^n)$ and
$$
\intr\A(x,\nabla u)\,dx\le L ;
$$
in particular, $u\in W^{1,p}(\R^n)$ if moreover $u\in L^p(\R^n)$;
\item[(iii)] $\E_s$ $\Gamma$-converges in $L^p(\R^n)$, as $s\uparrow1$, to the functional equal to $\intr\A(x,\nabla u)\,dx$ on $W^{1,p}(\R^n)$ and to $+\infty$ otherwise, and the constant sequence is a recovery sequence. The same holds for the restrictions to $\{u=0\text{ a.e.\ in }\R^n\setminus U\}$, for any open set $U$. Moreover, if $s_k\uparrow1$ and $\sup_k(\|u_k\|_p+\E_{s_k}(u_k))<\infty$, then $\{u_k\}$ is precompact in $L^p_{\rm loc}(\R^n)$ and its limit points belong to $W^{1,p}(\R^n)$.
\end{enumerate}
\end{thm}

Note that in (ii) no bound of $u_k$ in $L^p$ is required. By (iii) and the fundamental theorem of $\Gamma$-convergence, the minimizers of $\E_s(v)-\int_Ufv$ among functions vanishing outside a bounded open set $U$, with $f\in L^{p'}(U)$, converge in $L^p(U)$ to the minimizer of $\intr\A(x,\nabla v)-\int_Ufv$ among the functions $v\in W^{1,p}(\R^n)$ vanishing a.e.\ outside $U$; the details are as in the proof of Corollary \ref{cor.min} below (the uniform Poincar\'e inequality used there does not need any regularity of $U$).

\begin{rem}\label{rem.israel}
Theorem \ref{teo.fijo} contains \cite[Theorems 3.3, 3.4 and 5.1]{FBS-Israel} (up to the normalization, see Section \ref{sec.tools}): the hypothesis there is $m(x,r\omega)=a(x,\omega)+O(r)$ uniformly in $\omega$, for each $x$, which implies \eqref{K} by Example \ref{ej.K}(i) below. The proof of the lower bound in \cite{FBS-Israel} relies on \cite[Lemma 3.5]{FBS-Israel}, which states that if $g(t_2)\le c\,g(t_1)$ for $t_1<t_2$ and $h$ is decreasing, then $\int_0^1t^rgh\ge\frac{r+1}{2^{r+1}}\int_0^1t^rg\int_0^1t^rh$. For $c>2^{r+1}$ this is false (take $g=1$ on $(0,\alpha)$, $g=c$ on $(\alpha,1)$ and $h=\chi_{(0,\alpha)}$, with $\alpha$ small), and in the application $r+1=p(1-s)\to0$ and $c=m_+/m_-$; with a constant depending on $c$ the resulting lower bound is not sharp. The proof below does not use it, nor any monotonicity in $s$ or in the scale: compactness only needs the unweighted seminorm, through \eqref{coerc} and Lemma \ref{lema.transl}, and the sharp lower bound comes from Lemma \ref{lema.dual}.
\end{rem}

\subsection{\texorpdfstring{The condition \eqref{K} in logarithmic scale}{The condition (K) in logarithmic scale}}
For $s\in(0,1)$ let
\begin{equation}\label{Js}
J_s(h)=\frac{p(1-s)}{n\omega_n}\,\frac{\chi_{B_1}(h)}{|h|^{n-p(1-s)}},\qquad\text{so that}\qquad(1-s)\,|h|^{-n-sp}=\frac{n\omega_n}p\,J_s(h)\,|h|^{-p}\quad\text{in }B_1 .
\end{equation}
$J_s$ is radial, $\int J_s=1$ and $\int_{\{|h|>r\}}J_s=1-r^{p(1-s)}\to0$ as $s\uparrow1$, for every $r\in(0,1)$. 

\noindent In polar coordinates, $J_s(h)\,dh=\frac1{n\omega_n}\,d\Hn(\omega)\,d\mu_s(t)$, where $\mu_s$ is the probability measure $p(1-s)t^{p(1-s)-1}dt$ on $(0,1)$.

\begin{lema}\label{lema.Kstar}
Assume \eqref{H1} and \eqref{K}. Then, for every $R>0$,
\begin{equation}\label{Kstar}
e_s(R):=\int_{B_R}\int_{B_1}J_s(h)\,\Big|m(y,h)-a\Big(y,\frac h{|h|}\Big)\Big|\,dh\,dy\ \longrightarrow\ 0\qquad\text{as }s\uparrow1 .
\end{equation}
\end{lema}

\begin{proof}
Let $\Theta_R(t)=\int_{B_R}\frac1{n\omega_n}\int_{\Sn}|m(y,t\omega)-a(y,\omega)|\,d\Hn\,dy$, so that $e_s(R)=\int_0^1\Theta_R\,d\mu_s$, $\Theta_R\le\varpi_R(\delta)$ a.e.\ in $(0,\delta)$ and $\Theta_R\le m_+|B_R|$. Since $\mu_s([\delta,1))=1-\delta^{p(1-s)}$,
$$
e_s(R)\le\varpi_R(\delta)+m_+|B_R|\big(1-\delta^{p(1-s)}\big)
$$
for every $\delta\in(0,1)$. Let $s\uparrow1$ and then $\delta\to0$.
\end{proof}

Only \eqref{Kstar} is used in the proof of Theorem \ref{teo.fijo}; see Remark \ref{rem.ms}.

\subsection{The lower bound by nonlocal duality}

\begin{lema}\label{lema.dual}
Let $s_k\uparrow1$, $J_k=J_{s_k}$, and let $m_k\colon\R^n\times\R^n\to[m_-,m_+]$ be measurable and such that, for every $R>0$,
\begin{equation}\label{ek}
e_k(R):=\int_{B_R}\int_{B_1}J_k(h)\,\Big|m_k(x,h)-a\Big(x,\frac h{|h|}\Big)\Big|\,dh\,dx\ \longrightarrow\ 0\qquad\text{as }k\to\infty .
\end{equation}
Let $u_k\to u$ in $L^1_{\rm loc}(\R^n)$ and
$$
L_k=\intr\int_{B_1}m_k(x,h)\,J_k(h)\,\frac{|u_k(x)-u_k(x-h)|^p}{|h|^p}\,dh\,dx .
$$
If $L:=\liminf_kL_k<\infty$, then $u\in W^{1,1}_{\rm loc}(\R^n)$, $\nabla u\in L^p(\R^n;\R^n)$ and
$$
\intr\frac1{n\omega_n}\int_{\Sn}a(x,\omega)\,|\nabla u(x)\cdot\omega|^p\,d\Hn(\omega)\,dx\le L .
$$
\end{lema}

\begin{proof}
Let $D_k(x,h)=\frac{u_k(x)-u_k(x-h)}{|h|}$, and write $\omega=\frac h{|h|}$. Passing to a subsequence (which still satisfies \eqref{ek}) we may assume that $L_k\to L$ and $L_k\le L+1$. Since $m_k\ge m_-$,
\begin{equation}\label{cotaD}
\intr\int_{B_1}J_k|D_k|^p\,dh\,dx\le\frac{L+1}{m_-} .
\end{equation}

\emph{Step 1: Young's inequality.} For $t,\phi\in\R$,
\begin{equation}\label{young}
|t|^p\ge p\,\phi\,t-(p-1)|\phi|^{p'},\qquad\text{with equality if }\phi=|t|^{p-2}t .
\end{equation}
Let $R>0$ and let $\psi(x,\omega)=\sum_i\varphi_i(x)\beta_i(\omega)$ be a finite sum with $\varphi_i\in C^2_c(B_R)$ and $\beta_i$ bounded and measurable on $\Sn$; we denote by $\mathcal T_R$ the set of such functions. Applying \eqref{young} with $t=D_k$ and $\phi=\psi(x,\omega)/a(x,\omega)$, multiplying by $m_kJ_k$ and integrating over $B_R\times B_1$, since the integrand of $L_k$ is nonnegative,
$$
L_k\ge p\,I_k-p\,E_k-(p-1)\,I\!I_k,
$$
where all the integrals are over $B_R\times B_1$ and
$$
I_k=\iint J_k\,\psi\,D_k,\qquad E_k=\iint J_k\Big(1-\frac{m_k}a\Big)\psi\,D_k,\qquad I\!I_k=\iint J_k\,\frac{m_k}{a^{p'}}\,|\psi|^{p'} .
$$
They are finite by \eqref{cotaD} and H\"older's inequality.

\emph{Step 2: the coefficient.} This is the only place where the coefficient enters. Since $|1-\frac{m_k}a|\le\frac{|m_k-a|}{m_-}$ and $|m_k-a|^{p'}\le m_+^{p'-1}|m_k-a|$, H\"older's inequality with respect to $J_k\,dh\,dx$ and \eqref{cotaD} give
\begin{align*}
|E_k|&\le\frac{\|\psi\|_\infty}{m_-}\Big(\iint J_k|m_k-a|^{p'}\Big)^{1/p'}\Big(\iint J_k|D_k|^p\Big)^{1/p}\\
&\le\frac{\|\psi\|_\infty}{m_-}\big(m_+^{p'-1}e_k(R)\big)^{1/p'}\Big(\frac{L+1}{m_-}\Big)^{1/p}\to0 .
\end{align*}
Moreover $\frac{m_k}{a^{p'}}-a^{1-p'}=\frac{m_k-a}{a^{p'}}$, so $\big|I\!I_k-\iint J_k\,a^{1-p'}|\psi|^{p'}\big|\le m_-^{-p'}\|\psi\|_\infty^{p'}e_k(R)\to0$, and since $J_k$ is radial, supported in $B_1$ and $\int J_k=1$, the last double integral equals $\mathcal Q(\psi)$ for every $k$. Hence
$$
I\!I_k\ \longrightarrow\ \mathcal Q(\psi):=\intr\frac1{n\omega_n}\int_{\Sn}a^{1-p'}\,|\psi|^{p'}\,d\Hn\,dx .
$$
No continuity of $m_k$ or $a$ in $x$ has been used.

\emph{Step 3: discrete integration by parts.} Let $M_1=\sup|\nabla_x\psi|$ and $M_2=\sup\|D^2_x\psi\|$. For $0<|h|<1$, the change of variables $y=x-h$ gives
$$
\intr\psi(x,\omega)\big(u_k(x)-u_k(x-h)\big)\,dx=\intr u_k(y)\big(\psi(y,\omega)-\psi(y+h,\omega)\big)\,dy .
$$
Since $|\psi(y+h,\omega)-\psi(y,\omega)-\nabla_x\psi(y,\omega)\cdot h|\le\frac{M_2}2|h|^2$, and both sides vanish for $y\notin B_{R+1}$,
$$
\Big|\intr\psi(x,\omega)D_k(x,h)\,dx-g_k(\omega)\Big|\le\frac{M_2}2\,|h|\,\|u_k\|_{L^1(B_{R+1})},
$$
where
$$
g_k(\omega):=-\intr u_k(y)\,\nabla_x\psi(y,\omega)\cdot\omega\,dy .
$$
Multiplying by $J_k(h)$ and integrating over $B_1$,
$$
\Big|I_k-\frac1{n\omega_n}\int_{\Sn}g_k\,d\Hn\Big|\le\frac{M_2}2\,\sup_j\|u_j\|_{L^1(B_{R+1})}\int_{B_1}J_k(h)|h|\,dh ,
$$
and $\int_{B_1}J_k|h|\le r+\int_{\{|h|>r\}}J_k$ for every $r\in(0,1)$, so the right-hand side tends to $0$. Since $u_k\to u$ in $L^1(B_{R+1})$, $g_k\to g$ uniformly, where $g$ is defined as $g_k$ with $u$ in place of $u_k$. Hence $I_k\to\ell(\psi)$, with
$$
\ell(\psi):=-\intr\frac1{n\omega_n}\int_{\Sn}u(x)\,\nabla_x\psi(x,\omega)\cdot\omega\,d\Hn\,dx ,
$$
and we have proved that
\begin{equation}\label{dual}
L\ge p\,\ell(\psi)-(p-1)\,\mathcal Q(\psi)\qquad\text{for every }\psi\in\mathcal T:=\textstyle\bigcup_R\mathcal T_R .
\end{equation}

\emph{Step 4: $u\in W^{1,1}_{\rm loc}$ and $\nabla u\in L^p$.} Applying \eqref{dual} to $t\psi$, $t>0$, we get $L\ge p\,t\,\ell(\psi)-(p-1)t^{p'}\mathcal Q(\psi)$. If $\ell(\psi)>0$ then $\mathcal Q(\psi)>0$ (otherwise, letting $t\to\infty$ we would get $L=\infty$), and the choice $t=(\ell(\psi)/\mathcal Q(\psi))^{p-1}$ gives $L\ge\ell(\psi)^p\mathcal Q(\psi)^{1-p}$. Applying this also to $-\psi$,
$$
|\ell(\psi)|\le L^{1/p}\,\mathcal Q(\psi)^{1/p'} .
$$
Take $\psi(x,\omega)=\varphi(x)\omega_j$ with $\varphi\in C^\infty_c(\R^n)$. Since $\frac1{n\omega_n}\int_{\Sn}\omega_i\omega_j\,d\Hn=\frac1n\delta_{ij}$, we have $\ell(\psi)=-\frac1n\intr u\,\partial_j\varphi$, and $\mathcal Q(\psi)\le m_-^{1-p'}\|\varphi\|_{p'}^{p'}$. Hence
$$
\Big|\intr u\,\partial_j\varphi\,dx\Big|\le n\,m_-^{-1/p}\,L^{1/p}\,\|\varphi\|_{L^{p'}(\R^n)} ,
$$
and by the Riesz representation theorem the distributional derivative $\partial_ju$ belongs to $L^p(\R^n)$. Thus $u\in W^{1,1}_{\rm loc}(\R^n)$, $\nabla u\in L^p$, and
$$
\ell(\psi)=\intr\frac1{n\omega_n}\int_{\Sn}\psi\,\nabla u\cdot\omega\,d\Hn\,dx\qquad\text{for }\psi\in\mathcal T .
$$

\emph{Step 5: the supremum.} By \eqref{dual}, for every $\psi\in\mathcal T_R$,
$$
L\ge\Phi(\psi):=\intr\frac1{n\omega_n}\int_{\Sn}\Big(p\,\psi\,\nabla u\cdot\omega-(p-1)\,a^{1-p'}|\psi|^{p'}\Big)\,d\Hn\,dx .
$$
The functional $\Phi$ is continuous on $L^{p'}(B_R\times\Sn)$, since $\nabla u\in L^p$ and $a^{1-p'}\le m_-^{1-p'}$, and $\mathcal T_R$ is dense in that space (finite sums of products $\chi_E(x)\chi_F(\omega)$ are dense, and $\chi_E$ is a limit in $L^{p'}(B_R)$ of functions in $C^\infty_c(B_R)$). Hence $L\ge\Phi(\psi^*)$ for
$$
\psi^*=\chi_{B_R}\,a\,|\nabla u\cdot\omega|^{p-2}\,\nabla u\cdot\omega\in L^{p'}(B_R\times\Sn).
$$
By the equality case in \eqref{young} with $t=\nabla u\cdot\omega$ and $\phi=\psi^*/a$, the integrand of $\Phi(\psi^*)$ is $a|\nabla u\cdot\omega|^p$ in $B_R$. Therefore $L\ge\int_{B_R}\frac1{n\omega_n}\int_{\Sn}a|\nabla u\cdot\omega|^p$, and the result follows letting $R\to\infty$.
\end{proof}

\subsection{Proof of Theorem \ref{teo.fijo}}

\begin{proof}[Proof of Theorem \ref{teo.fijo}]
(i) Let first $v\in C^2_c(\R^n)$ with $\operatorname{supp}v\subset B_R$, $R\ge1$, and let $M_1=\|\nabla v\|_\infty$, $M_2=\|D^2v\|_\infty$. By Lemma \ref{lema.cota}, $\E_s(v;1,\infty)\le C(1-s)\|v\|_p^p\to0$. For $|h|<1$ only the points $x\in B_{R+1}$ contribute to $\E_s(v;0,1)$. Let $\alpha=v(x)-v(x-h)$ and $\beta=\nabla v(x)\cdot h$; then $|\alpha|,|\beta|\le M_1|h|$ and $|\alpha-\beta|\le\frac{M_2}2|h|^2$, so $\big||\alpha|^p-|\beta|^p\big|\le\frac p2M_1^{p-1}M_2|h|^{p+1}$. Hence replacing $|\alpha|^p$ by $|\beta|^p$ in $\E_s(v;0,1)$ costs at most
$$
(1-s)\,m_+\,\frac p2\,M_1^{p-1}M_2\,|B_{R+1}|\int_{B_1}\frac{dh}{|h|^{n-1-p(1-s)}}\le C_v\,(1-s) ,
$$
and then replacing $m(x,h)$ by $a(x,\frac h{|h|})$ costs at most
$$
(1-s)\,M_1^p\int_{B_{R+1}}\int_{B_1}\Big|m(x,h)-a\Big(x,\frac h{|h|}\Big)\Big|\frac{dh\,dx}{|h|^{n-p(1-s)}}=\frac{n\omega_n}p\,M_1^p\,e_s(R+1)\ \longrightarrow\ 0
$$
by Lemma \ref{lema.Kstar}. Finally, in polar coordinates, since $(1-s)\int_0^1r^{p(1-s)-1}dr=\frac1p$,
$$
(1-s)\intr\int_{B_1}a\Big(x,\frac h{|h|}\Big)\frac{|\nabla v(x)\cdot h|^p}{|h|^{n+sp}}\,dh\,dx=\intr\A(x,\nabla v)\,dx .
$$
So $\E_s(v)\to\intr\A(x,\nabla v)$. For $u\in W^{1,p}(\R^n)$ the result follows by density, since $\E_s^{1/p}$ is a seminorm bounded by $C\|\cdot\|_{W^{1,p}}$ uniformly in $s\ge\frac12$ (Lemma \ref{lema.cota}), and $v\mapsto(\intr\A(x,\nabla v))^{1/p}$ is a seminorm bounded by $c_2^{1/p}\|\nabla v\|_p$.

(ii) Pass to a subsequence realizing the liminf. By \eqref{Js}, keeping only $|h|<1$,
$$
\E_{s_k}(u_k)\ge\frac{n\omega_n}p\intr\int_{B_1}m(x,h)\,J_{s_k}(h)\,\frac{|u_k(x)-u_k(x-h)|^p}{|h|^p}\,dh\,dx ,
$$
and by Lemma \ref{lema.Kstar} the hypothesis \eqref{ek} holds with $m_k=m$. Lemma \ref{lema.dual} gives $u\in W^{1,1}_{\rm loc}$, $\nabla u\in L^p$ and
$$
L\ge\frac{n\omega_n}p\intr\frac1{n\omega_n}\int_{\Sn}a\,|\nabla u\cdot\omega|^p\,d\Hn\,dx=\intr\A(x,\nabla u)\,dx .
$$

(iii) The liminf inequality is (ii), and the constant sequence is a recovery sequence by (i). If the $u_k$ vanish outside $U$ and $u_k\to u$ in $L^p$, then $u$ vanishes outside $U$, and if $u\in W^{1,p}(\R^n)$ vanishes outside $U$ the constant sequence is admissible; no regularity of $U$ is needed. For the compactness, let $\Lambda_0=\sup_k\E_{s_k}(u_k)$. By \eqref{coerc} and Lemma \ref{lema.transl}, $\|\tau_hu_k-u_k\|_p^p\le C_{n,p}|h|^{s_kp}\Lambda_0/m_-$ for every $h$ and $k$. Since $s_k\ge\frac12$ for $k$ large, the sequence is bounded and equicontinuous in $L^p$ with respect to translations, and the Riesz--Fr\'echet--Kolmogorov theorem gives the precompactness in $L^p_{\rm loc}(\R^n)$. The limit points belong to $L^p(\R^n)$ by Fatou's lemma, hence to $W^{1,p}(\R^n)$ by (ii).
\end{proof}

\subsection{Examples}

\begin{exam}\label{ej.K}
In all the examples $m$ takes values in $[m_-,m_+]$.
\begin{enumerate}
\item[(i)] \emph{Pointwise convergence.} Assume that for a.e.\ $y$,
$$
\eta(y,r):=\operatorname*{ess\,sup}_{0<|h|<r}\Big|m(y,h)-a\Big(y,\frac h{|h|}\Big)\Big|\ \longrightarrow\ 0\qquad\text{as }r\to0 .
$$
This is the case if $m(y,r\omega)=a(y,\omega)+O(r)$ uniformly in $\omega$ for each $y$, as in \cite{FBS-Israel}, and in particular under \eqref{H3}. Then \eqref{K} holds, since $\varpi_R(r)\le\int_{B_R}\eta(y,r)\,dy\to0$ by dominated convergence (the function $\eta(\cdot,r)$ is measurable, being the limit as $q\to\infty$ of the $L^q$ averages of $h\mapsto|m(y,h)-a(y,\frac h{|h|})|$ over $B_r$). No uniformity in $y$ and no rate is needed.
\item[(ii)] \emph{Coefficients built from a measurable function.} Let $c\colon\R^n\to[\alpha,\beta]$ be measurable, $G\colon[\alpha,\beta]^2\to[m_-,m_+]$ continuous, and
$$
m(x,h)=G\big(c(x),c(x-h)\big),\qquad a(x,\omega)=G\big(c(x),c(x)\big).
$$
For instance $G(t_1,t_2)=\frac12(t_1+t_2)$, $\sqrt{t_1t_2}$ or $\max\{t_1,t_2\}$, which give symmetric kernels. By uniform continuity, given $\nu>0$ there is $\gamma>0$ with $|G(t_1,t_2)-G(t_1,t_1)|\le\nu+\frac{m_+}\gamma|t_2-t_1|$, so
$$
\varpi_R(r)\le\nu\,|B_R|+\frac{m_+}\gamma\sup_{|z|\le r}\int_{B_R}|c(x+z)-c(x)|\,dx ,
$$
and the last term tends to $0$ as $r\to0$ by the continuity of translations in $L^1(\R^n)$, applied to $c\,\chi_{B_{R+1}}$. Hence \eqref{K} holds. On the other hand, the condition of (i), and hence \eqref{H3}, fails in general: let $c=1+\chi_F$, where $F$ is a measurable set with $0<|F\cap B|<|B|$ for every ball $B$, and assume $G(1,2)\ne G(1,1)$. For every $y\notin F$ and every $r>0$ the set $\{h\in B_r\colon c(y-h)=2\}$ has positive measure, so $\eta(y,r)=|G(1,2)-G(1,1)|>0$ $\forall r$. Since $|\R^n\setminus F|>0$, the condition of (i) fails on a set of positive measure.
\item[(iii)] \emph{Midpoint coefficients.} Let $b\colon\R^n\times\Sn\to[m_-,m_+]$ be measurable and
$$
m(x,h)=b\Big(x-\frac h2,\frac h{|h|}\Big),\qquad a(x,\omega)=b(x,\omega).
$$
For a.e.\ $\omega$, $\sup_{|z|\le r}\int_{B_R}|b(x+z,\omega)-b(x,\omega)|\,dx\to0$ as $r\to0$, by the continuity of translations in $L^1$, and \eqref{K} follows by dominated convergence in $\omega$.
\end{enumerate}
\end{exam}

\begin{rem}[Kernels depending on $s$]\label{rem.ms}
The proof of Theorem \ref{teo.fijo} only uses \eqref{H1} and \eqref{Kstar}. Hence the theorem holds, with the same proof, for the energies of a family of kernels $m_s$ satisfying \eqref{H1} with constants independent of $s$ and
\begin{equation}\label{Kstar2}
\lim_{s\uparrow1}\int_{B_R}\int_{B_1}J_s(h)\,\Big|m_s(y,h)-a\Big(y,\frac h{|h|}\Big)\Big|\,dh\,dy=0\qquad\text{for every }R>0 .
\end{equation}
This holds, for instance, if $m_s\to m$ uniformly in $\{y\in B_R,\ 0<|h|<\delta_R\}$ for every $R$, where $m$ satisfies \eqref{K}. Condition \eqref{Kstar2} asks that $m_s(y,t\omega)$ be close to $a(y,\omega)$, in mean, for most scales $t$ with respect to the measure $\mu_s$; since each dyadic interval has $\mu_s$-measure at most $p(1-s)\log2$, a discrepancy confined to $o(\frac1{1-s})$ dyadic scales is invisible.

The convergence $m_s\to m$ a.e., or weakly-$*$ in $L^\infty$, is not enough. Let $m_s(y,h)=\alpha_s(|h|)$ with $\alpha_s\colon(0,\infty)\to[m_-,m_+]$ measurable. For $u\in W^{1,p}(\R^n)$ the functions $\phi(t,\omega)=t^{-p}\|\tau_{-t\omega}u-u\|_p^p$ satisfy $\phi(t,\omega)\le\|\nabla u\cdot\omega\|_p^p$ and $\phi(t,\omega)\to\|\nabla u\cdot\omega\|_p^p$ as $t\to0$, uniformly in $\omega$, and one obtains, in polar coordinates,
$$
\lim_{s\uparrow1}\Big(\E_s(u)-\frac{K_{n,p}}p\,\|\nabla u\|_p^p\int_0^1\alpha_s\,d\mu_s\Big)=0 .
$$
If $\alpha_s=m_-+(m_+-m_-)\chi_{(0,r_s)}$ with $r_s=e^{-c/(p(1-s))}$, $c>0$, then $\int_0^1\alpha_s\,d\mu_s=m_-+(m_+-m_-)e^{-c}$, although $m_s\to m_-$ a.e. Here $r_s$ plays the role of $\ve$ and $r_s^{p(1-s)}=e^{-c}$ that of $\theta$ in Theorem \ref{teo.main}: this is the same mechanism, in a translation invariant situation. Finally, neither \eqref{Kstar2} nor \eqref{K} is necessary for the convergence of $\E_s(u)$, $u\in W^{1,p}(\R^n)$: if $\alpha=m_-+(m_+-m_-)\chi_F$ with $F=\bigcup_{j\text{ even}}(2^{-j-1},2^{-j}]$, then $\mu_s(F)\to\frac12$ and $\E_s(u)\to\frac{m_-+m_+}2\frac{K_{n,p}}p\|\nabla u\|_p^p$, while \eqref{Kstar2} fails for every profile.
\end{rem}

\section{The joint limit}\label{sec.conjunto}

\subsection{Setting and statements}
In this section $m$ is Borel measurable, $Q$-periodic in its first variable, and satisfies \eqref{H1} and \eqref{H3}. Then $a$ is Borel measurable (being the pointwise limit of $m(y,r_j\omega)$ as $r_j\to0$), $Q$-periodic in $y$, and $m_-\le a\le m_+$. The Borel measurability of $m$ is needed here: for fixed $r$, the function $(y,\omega)\mapsto m(y,r\omega)$ is the restriction of $m$ to the set $\{|h|=r\}$, which is Lebesgue null in $\R^n\times\R^n$, so it need not be measurable if $m$ is only Lebesgue measurable. We write
$$
\bar m(h)=\int_Qm(y,h)\,dy,\qquad\bar a(\omega)=\int_Qa(y,\omega)\,dy,
$$
and $\bar\E_s(u;\rho_1,\rho_2)$ for the energy defined as $\E_{\ve,s}(u;\rho_1,\rho_2)$ with $m(\frac x\ve,h)$ replaced by $\bar m(h)$. Note that \eqref{H3} gives $|\bar m(r\omega)-\bar a(\omega)|\le C_0r$ for $0<r<1$. Besides $\A$ we consider the densities
\begin{align*}
\bar\A(\xi)&=\int_Q\A(y,\xi)\,dy=\frac1p\int_{\Sn}\bar a(\omega)|\xi\cdot\omega|^p\,d\Hn,\\
\A^*(\xi)&=\inf\Big\{\int_Q\A(y,\xi+\nabla v(y))\,dy\colon v\in W^{1,p}_{\rm per}(Q)\Big\},
\end{align*}
which satisfy
\begin{equation}\label{growth}
c_1|\xi|^p\le\A^*(\xi)\le\bar\A(\xi)\le c_2|\xi|^p .
\end{equation}
We consider sequences
\begin{equation}\label{theta}
\ve_k\to0,\qquad s_k\to1,\qquad\ve_k^{p(1-s_k)}\to\theta\in[0,1],
\end{equation}
that is, $(1-s_k)\log(1/\ve_k)\to\lambda\in[0,\infty]$ and $\theta=e^{-p\lambda}$. We write $\E_k=\E_{\ve_k,s_k}$, and
$$
\A_\theta=\theta\,\A^*+(1-\theta)\,\bar\A,\qquad\E_\theta(u)=\begin{cases}\displaystyle\intr\A_\theta(\nabla u)\,dx&u\in W^{1,p}(\R^n),\\+\infty&\text{otherwise.}\end{cases}
$$

\begin{thm}\label{teo.main}
Let $m$ be Borel measurable and $Q$-periodic in its first variable, and assume \eqref{H1}, \eqref{H3} and \eqref{theta}. Then:
\begin{enumerate}
\item (Compactness) if $\{u_k\}$ is bounded in $L^p(\R^n)$ and $\sup_k\E_k(u_k)<\infty$, then $\{u_k\}$ is precompact in $L^p_{\rm loc}(\R^n)$ and every limit point belongs to $W^{1,p}(\R^n)$;
\item (Liminf) if $\{u_k\}$ is bounded in $L^p(\R^n)$ and $u_k\to u$ in $L^p_{\rm loc}(\R^n)$, then
$$
\E_\theta(u)\le\liminf_{k\to\infty}\E_k(u_k);
$$
\item (Limsup) for every $u\in W^{1,p}(\R^n)$ there exist $u_k\to u$ in $L^p(\R^n)$ with
$$
\limsup_{k\to\infty}\E_k(u_k)\le\E_\theta(u).
$$
If $u$ vanishes a.e.\ outside a bounded Lipschitz open set $\Omega$, the $u_k$ can be taken with compact support in $\Omega$.
\end{enumerate}
In particular, $\E_k$ $\Gamma$-converges to $\E_\theta$ in $L^p(\R^n)$, and the same holds for the functionals restricted to $\{u=0\text{ a.e.\ in }\R^n\setminus\Omega\}$ when $\Omega$ is bounded and Lipschitz.
\end{thm}

The lower bound splits into two independent statements, one for each range of scales. Given $s_k\to1$ we fix once and for all
\begin{equation}\label{Rk}
R_k=\exp\big((1-s_k)^{-1/2}\big),\qquad\text{so that}\qquad R_k\to\infty\quad\text{and}\quad R_k^{\pm p(1-s_k)}\to1 .
\end{equation}
Under \eqref{theta}, both $(R_k\ve_k)^{p(1-s_k)}$ and $(\ve_k/R_k)^{p(1-s_k)}$ converge to $\theta$.

\begin{thm}[Large scales]\label{teo.macro}
Assume \eqref{H1}, \eqref{H3} and \eqref{theta} with $\theta<1$. Let $\{u_k\}$ be bounded in $L^p(\R^n)$ with $u_k\to u$ in $L^p_{\rm loc}(\R^n)$ and $\sup_k\E_k(u_k)<\infty$. Then $u\in W^{1,p}(\R^n)$ and
$$
(1-\theta)\intr\bar\A(\nabla u)\,dx\le\liminf_{k\to\infty}\E_k(u_k;R_k\ve_k,\infty).
$$
\end{thm}

\begin{thm}[Small scales]\label{teo.micro}
Assume \eqref{H1}, \eqref{H3} and \eqref{theta} with $\theta>0$. Let $\{u_k\}$ be bounded in $L^p(\R^n)$ with $u_k\to u$ in $L^p_{\rm loc}(\R^n)$ and $\sup_k\E_k(u_k)<\infty$. Then $u\in W^{1,p}(\R^n)$ and
$$
\theta\intr\A^*(\nabla u)\,dx\le\liminf_{k\to\infty}\E_k(u_k;0,\ve_k/R_k).
$$
\end{thm}

The proof of Theorems \ref{teo.macro} and \ref{teo.micro} will be delivered in Sections \ref{sec.macro} and \ref{sec.micro}, respectively.

\begin{cor}\label{cor.min}
Assume \eqref{H1}, \eqref{H3} and \eqref{theta}. Let $\Omega$ be a bounded Lipschitz open set, $f\in L^{p'}(\Omega)$, and let $u_k$ be the minimizer of $\E_k(v)-\int_\Omega fv$ among $v\in L^p(\R^n)$ with $v=0$ a.e.\ in $\R^n\setminus\Omega$. Then $u_k\to u_\theta$ in $L^p(\Omega)$, where $u_\theta\in W^{1,p}_0(\Omega)$ is the unique minimizer of $\int_\Omega\A_\theta(\nabla v)-\int_\Omega fv$ in $W^{1,p}_0(\Omega)$.
\end{cor}

\begin{rem}
The density $\A^*$ is of class $C^1$: it is the marginal $\A^*(\xi)=\min_w\int_Q\A(y,\xi+\nabla w)$ of a jointly convex, G\^ateaux differentiable functional, so its subdifferential at each point is the singleton $\{\bar\tau\}$ of Lemma \ref{lema.celda}. Hence $u_\theta$ is the weak solution of $-\operatorname{div}(\nabla_\xi\A_\theta(\nabla u))=f$ in $\Omega$, $u=0$ on $\partial\Omega$.
\end{rem}

\begin{rem}
If $\ve_k^{p(1-s_k)}$ has two different limit points $\theta\ne\theta'$ and $\A^*\ne\bar\A$, the corresponding subsequences have different $\Gamma$-limits, so $\E_k$ does not $\Gamma$-converge (recall that $\A^*\le\bar\A$, so $\A^*\ne\bar\A$ gives $\E_\theta\ne\E_{\theta'}$). The two iterated limits correspond to $\theta=1$ ($s\uparrow1$ first) and $\theta=0$ ($\ve\downarrow0$ first).
\end{rem}

In the quadratic case, Theorem \ref{teo.main} applies to the nonlocal approximation of the energies of general periodic matrices.

\begin{cor}[Periodic matrices]\label{cor.matrices}
Let $p=2$ and let $A\colon\R^n\to\R^{n\times n}$ be Borel measurable, $Q$-periodic, with values in the symmetric matrices, and such that $\alpha|\xi|^2\le A(y)\xi\cdot\xi\le\beta|\xi|^2$ for some $0<\alpha\le\beta$. Let
$$
m(y,h)=\frac{2}{\omega_n\sqrt{\det A(y)}}\Big(A(y)^{-1}\frac h{|h|}\cdot\frac h{|h|}\Big)^{-\frac{n+2}2},\qquad h\ne0.
$$
Then $\A(y,\xi)=A(y)\xi\cdot\xi$, and under \eqref{theta} the energies $\E_k$ $\Gamma$-converge in $L^2(\R^n)$ to the functional equal to $\intr A_\theta\nabla u\cdot\nabla u\,dx$ on $H^1(\R^n)$ and to $+\infty$ otherwise, where
$$
A_\theta=\theta\,A_{\rm hom}+(1-\theta)\int_QA(y)\,dy
$$
and $A_{\rm hom}$ is the homogenized matrix of $A$, that is, $A_{\rm hom}\xi\cdot\xi=\inf\{\int_QA(y)(\xi+\nabla v)\cdot(\xi+\nabla v)\,dy\colon v\in H^1_{\rm per}(Q)\}$.
\end{cor}

\begin{proof}
The kernel $m$ is Borel measurable, $Q$-periodic in $y$ and homogeneous of degree $0$ in $h$, so \eqref{H3} holds with $a(y,\omega)=m(y,\omega)$ and $C_0=0$, and \eqref{H1} holds with $m_-=\frac{2\alpha^{(n+2)/2}}{\omega_n\beta^{n/2}}$ and $m_+=\frac{2\beta^{(n+2)/2}}{\omega_n\alpha^{n/2}}$, because the eigenvalues of $A(y)$ belong to $[\alpha,\beta]$. (Recall that $K_{n,2}=\int_{\Sn}\omega_1^2\,d\Hn=\omega_n$.) To compute $\A$, fix $y$, write $A=A(y)$, and let $E=\{z\colon A^{-1}z\cdot z\le1\}=A^{1/2}(B_1)$, whose radial function is $r_E(\omega)=(A^{-1}\omega\cdot\omega)^{-1/2}$. The change of variables $z=A^{1/2}w$ and polar coordinates give, respectively,
\begin{align*}
\int_E(\xi\cdot z)^2dz&=\sqrt{\det A}\int_{B_1}(A^{1/2}\xi\cdot w)^2dw=\frac{\omega_n\sqrt{\det A}}{n+2}\,A\xi\cdot\xi,\\
\int_E(\xi\cdot z)^2dz&=\frac1{n+2}\int_{\Sn}r_E(\omega)^{n+2}(\xi\cdot\omega)^2d\Hn ,
\end{align*}
where we used that $\int_{B_1}(\eta\cdot w)^2dw=\frac{\omega_n}{n+2}|\eta|^2$. Since $a(y,\omega)=\frac{2}{\omega_n\sqrt{\det A}}r_E(\omega)^{n+2}$, we obtain $\A(y,\xi)=\frac12\int_{\Sn}a(y,\omega)(\xi\cdot\omega)^2d\Hn=A\xi\cdot\xi$. Hence $\bar\A(\xi)=(\int_QA)\xi\cdot\xi$ and $\A^*(\xi)=A_{\rm hom}\xi\cdot\xi$, and the result follows from Theorem \ref{teo.main}.
\end{proof}

\subsection{The lower bound from the large scales}\label{sec.macro}

\begin{lema}[Averaging]\label{lema.prom}
There is a constant $c_n$ such that for every $u\in L^p(\R^n)$, $\ve,\delta>0$ and $h\in\R^n$,
$$
\Big|\intr\Big(m\Big(\frac x\ve,h\Big)-\bar m(h)\Big)|u_\delta(x)-u_\delta(x-h)|^p\,dx\Big|\le c_n\,p\,c_\varrho\,m_+\,\frac\ve\delta\,\|\tau_{-h}u-u\|_p^p .
$$
Moreover, for every $v\in W^{1,p}(\R^n)$,
$$
\Big|\intr\Big(m\Big(\frac x\ve,h\Big)-\bar m(h)\Big)|v(x)-v(x-h)|^p\,dx\Big|\le c_n\,p\,m_+\,\ve\,\|\tau_{-h}v-v\|_p^{p-1}\|\tau_{-h}\nabla v-\nabla v\|_p .
$$
\end{lema}

\begin{proof}
Let $v\in W^{1,p}(\R^n)$ and $F(x)=|v(x)-v(x-h)|^p$. Then $F\in W^{1,1}(\R^n)$ and, by H\"older's inequality,
$$
\|\nabla F\|_1\le p\,\|\tau_{-h}v-v\|_p^{p-1}\|\tau_{-h}\nabla v-\nabla v\|_p .
$$
Let $g(y)=m(y,h)-\bar m(h)$, which is $Q$-periodic, has zero mean on $Q$ and satisfies $|g|\le m_+$. On each cube $Q_j=\ve(j+Q)$, $j\in\mathbb Z^n$, we have $\int_{Q_j}g(x/\ve)\,dx=0$, so by the Poincar\'e inequality in $L^1$ on cubes,
$$
\Big|\int_{Q_j}g\Big(\frac x\ve\Big)F\,dx\Big|=\Big|\int_{Q_j}g\Big(\frac x\ve\Big)(F-F_{Q_j})\,dx\Big|\le m_+\,c_n\,\ve\int_{Q_j}|\nabla F|\,dx ,
$$
where $F_{Q_j}$ is the average of $F$ on $Q_j$. Summing over $j$ gives the second inequality. The first one follows applying it to $v=u_\delta$ and using Lemma \ref{lema.moll}(2).
\end{proof}

\begin{lema}[BBM lower bound on a band, translation invariant kernels]\label{lema.banda}
Let $s_k\uparrow1$ and $\rho_k>0$ with $\rho_k^{p(1-s_k)}\to\vartheta\in[0,1)$. Let $\{v_k\}$ be bounded in $L^p(\R^n)$ with $v_k\to u$ in $L^p_{\rm loc}(\R^n)$. If $\liminf_k\bar\E_{s_k}(v_k;\rho_k,\infty)<\infty$, then $u\in W^{1,p}(\R^n)$ and
$$
(1-\vartheta)\intr\bar\A(\nabla u)\,dx\le\liminf_{k\to\infty}\bar\E_{s_k}(v_k;\rho_k,\infty).
$$
\end{lema}

\begin{proof}
Since $\vartheta<1$ and $s_k\uparrow1$, we have $\rho_k\to0$. Fix $\sigma>0$, $\kappa\in(0,1)$ and $M>0$, and let $w_k=(v_k)_\sigma$. Since $\bar m$ does not depend on $x$, Jensen's inequality and translation invariance give, for every $h$, $\intr|w_k(x)-w_k(x-h)|^pdx\le\intr|v_k(x)-v_k(x-h)|^pdx$, and therefore
$$
\bar\E_{s_k}(v_k;\rho_k,\infty)\ge\bar\E_{s_k}(w_k;\rho_k,\kappa).
$$
The sequence $\{w_k\}$ is bounded in $C^2(\R^n)$, say $\|\nabla w_k\|_\infty\le K_1$ and $\|D^2w_k\|_\infty\le K_2$, and $w_k\to u_\sigma$ in $C^1(\overline B_M)$. For $x\in\R^n$ and $h=r\omega$ with $0<r<\kappa$, the numbers $\alpha=\frac{|w_k(x)-w_k(x-h)|}{r}$ and $\beta=|\nabla w_k(x)\cdot\omega|$ belong to $[0,K_1]$ and $|\alpha-\beta|\le\frac{K_2}2r$, so $\alpha^p\ge\beta^p-\frac p2K_1^{p-1}K_2\,r$. Using also $\bar m(r\omega)\ge\bar a(\omega)-C_0r$ and $\bar m\le m_+$, for $k$ large (so that $\rho_k<\kappa$),
\begin{align*}
(1-s_k)\int_{\{\rho_k<|h|<\kappa\}}&\bar m(h)\frac{|w_k(x)-w_k(x-h)|^p}{|h|^{n+s_kp}}\,dh\\
&=(1-s_k)\int_{\rho_k}^\kappa\int_{\Sn}\bar m(r\omega)\,\alpha^p\,r^{p(1-s_k)-1}\,d\Hn\,dr\\
&\ge\bar\A(\nabla w_k(x))\,\big(\kappa^{p(1-s_k)}-\rho_k^{p(1-s_k)}\big)-C(1-s_k),
\end{align*}
where $C$ depends only on $n,p,m_+,C_0,K_1,K_2$, and we used $\int_0^\kappa r^{p(1-s_k)}dr\le1$ and
$$
(1-s_k)\int_{\rho_k}^\kappa r^{p(1-s_k)-1}dr=\frac1p\big(\kappa^{p(1-s_k)}-\rho_k^{p(1-s_k)}\big).
$$
Integrating over $x\in B_M$ and letting $k\to\infty$, since $\kappa^{p(1-s_k)}\to1$ and $\nabla w_k\to\nabla u_\sigma$ uniformly on $B_M$,
$$
L:=\liminf_{k\to\infty}\bar\E_{s_k}(v_k;\rho_k,\infty)\ge(1-\vartheta)\int_{B_M}\bar\A(\nabla u_\sigma)\,dx .
$$
Letting $M\to\infty$ we get $(1-\vartheta)\intr\bar\A(\nabla u_\sigma)\le L$ for every $\sigma>0$. By \eqref{growth}, $\{\nabla u_\sigma\}_\sigma$ is bounded in $L^p$; since $u\in L^p(\R^n)$ by Fatou's lemma and $u_\sigma\to u$ in $L^p$, we conclude that $u\in W^{1,p}(\R^n)$, and then $\nabla u_\sigma=(\nabla u)_\sigma\to\nabla u$ in $L^p$. Since $\bar\A$ is continuous and has $p$-growth, $\intr\bar\A(\nabla u_\sigma)\to\intr\bar\A(\nabla u)$ as $\sigma\to0$.
\end{proof}

\begin{proof}[Proof of Theorem \ref{teo.macro}]
Let $\Lambda_0=\sup_k\E_k(u_k)$, $\rho_k=R_k\ve_k$ and $\delta_k=\sqrt{\ve_k\rho_k}=\sqrt{R_k}\,\ve_k$, so that
$$
\frac{\ve_k}{\delta_k}=\frac{\delta_k}{\rho_k}=R_k^{-1/2}\to0 .
$$
Since $\theta<1$, $\rho_k\to0$ and hence $\delta_k\to0$. Let $v_k=(u_k)_{\delta_k}$. By \eqref{coerc}, $(1-s_k)[u_k]_{s_k,p}^p\le\Lambda_0/m_-$. We may assume $s_k\ge\frac12$.

\emph{Step 1: from $u_k$ to $v_k$.} The quantity $\E_k(\cdot\,;\rho_k,\infty)^{1/p}$ is a seminorm (it is an $L^p$ norm of $(x,h)\mapsto u(x)-u(x-h)$ with respect to a positive measure). By \eqref{H1} and Lemma \ref{lema.moll}(3) with $\delta_k\le\rho_k$,
$$
\E_k(u_k-v_k;\rho_k,\infty)\le m_+\,\frac{2^pn\omega_nC_{n,p}}{s_kp}(1-s_k)\frac{\Lambda_0}{m_-}\le C(1-s_k)\Lambda_0 ,
$$
and therefore
\begin{equation}\label{paso1}
\E_k(u_k;\rho_k,\infty)^{1/p}\ge\E_k(v_k;\rho_k,\infty)^{1/p}-\big(C(1-s_k)\Lambda_0\big)^{1/p}.
\end{equation}

\emph{Step 2: from $m(x/\ve_k,h)$ to $\bar m(h)$.} Multiplying the first inequality in Lemma \ref{lema.prom} by $(1-s_k)|h|^{-n-s_kp}$, integrating over $\{|h|>\rho_k\}$ and using \eqref{H1},
\begin{equation}\label{paso2}
\big|\E_k(v_k;\varrho_k,\infty)-\bar\E_{s_k}(v_k;\varrho_k,\infty)\big|\le c_n\,p\,c_\varrho\,\frac{m_+}{m_-}\,\frac{\ve_k}{\delta_k}\,\E_k(u_k;\varrho_k,\infty)\le C\,R_k^{-1/2}\Lambda_0 .
\end{equation}

\emph{Step 3: conclusion.} By Lemma \ref{lema.moll}(1), $\|v_k-u_k\|_p^p\le C\delta_k^{s_kp}\Lambda_0\to0$, so $\{v_k\}$ is bounded in $L^p$ and $v_k\to u$ in $L^p_{\rm loc}$. By \eqref{paso1} and \eqref{paso2},
$$
\liminf_{k\to\infty}\E_k(u_k;\rho_k,\infty)\ge\liminf_{k\to\infty}\bar\E_{s_k}(v_k;\rho_k,\infty),
$$
and the right-hand side is finite. Lemma \ref{lema.banda} with $\vartheta=\theta$ concludes the proof.
\end{proof}

\subsection{The lower bound from the small scales}\label{sec.micro}

The proof has two ingredients: a reduction of the small scales to a local energy with a regularized weight (Proposition \ref{prop.red}), and the liminf inequality of periodic homogenization for local functionals whose coefficients also vary along the sequence (Proposition \ref{prop.homog}). The proof of the latter is the classical one by duality, which only uses weak convergence.

For $\eta>0$ let
\begin{equation}\label{beta}
b_\eta(y,\omega)=\Big(\int\varrho_\eta(z)\,a(y-z,\omega)^{-\frac1{p-1}}dz\Big)^{-(p-1)},\qquad\mathcal B_\eta(y,\xi)=\frac1p\int_{\Sn}b_\eta(y,\omega)|\xi\cdot\omega|^pd\Hn .
\end{equation}
Then $b_\eta$ is measurable, $m_-\le b_\eta\le m_+$, $b_\eta$ is $Q$-periodic in $y$ if $a$ is, and, by the Lebesgue differentiation theorem applied to $a(\cdot,\omega)^{-1/(p-1)}$ for each $\omega$ (recall that $a$ is Borel) and Fubini's theorem, $b_{\eta_k}\to a$ a.e.\ in $\R^n\times\Sn$ for every sequence $\eta_k\to0$.

\begin{prop}\label{prop.red}
Assume \eqref{H1} and \eqref{H3}. Let $\ve>0$, $s\in(0,1)$, $0<\delta<1$ with $C_0\delta<m_-$, $u\in L^p(\R^n)$, $v=u_\delta$ and $\Lambda=(1-s)[u]_{s,p}^p$. Then
$$
\E_{\ve,s}(u;0,\delta)^{1/p}\ge\Big(1-\frac{C_0\delta}{m_-}\Big)^{1/p}\Big[\Big(\delta^{p(1-s)}\intr\mathcal B_{\delta/\ve}\Big(\frac x\ve,\nabla v\Big)dx\Big)^{1/p}-\big(C_1(1-s)\Lambda\big)^{1/p}\Big],
$$
where $C_1$ depends only on $n$, $p$, $m_+$ and $\varrho$.
\end{prop}

\begin{proof}
For a weight $b(y,\omega)$ we write
$$
\E^b(w)=(1-s)\int_{\{|h|<\delta\}}\intr b\Big(\frac x\ve,\frac h{|h|}\Big)\frac{|w(x)-w(x-h)|^p}{|h|^{n+sp}}\,dx\,dh .
$$

\emph{Step 1: from $m$ to $a$.} For $|h|<\delta<1$, \eqref{H3} and $a\ge m_-$ give $m(y,h)\ge a(y,\frac h{|h|})-C_0\delta\ge(1-\frac{C_0\delta}{m_-})a(y,\frac h{|h|})$. Hence $\E_{\ve,s}(u;0,\delta)\ge(1-\frac{C_0\delta}{m_-})\E^a(u)$.

\emph{Step 2: mollification.} Let $\eta=\delta/\ve$. For fixed $\omega$, the function $b(x)=a(\frac x\ve,\omega)$ satisfies $b^\delta(x)=b_\eta(\frac x\ve,\omega)$, with the notation of Lemma \ref{lema.holder} (change variables $z\mapsto\ve z$). Applying that lemma for each $h$ and integrating, $\E^a(u)\ge\E^{b_\eta}(v)$.

\emph{Step 3: Taylor expansion in $L^p$.} Consider on $(0,\delta)\times\Sn\times\R^n$ the measure $d\mu=(1-s)\,b_\eta(\frac x\ve,\omega)\,r^{p(1-s)-1}\,dr\,d\Hn(\omega)\,dx$ and the functions
$$
F(r,\omega,x)=\frac{v(x)-v(x-r\omega)}r,\qquad G(r,\omega,x)=\nabla v(x)\cdot\omega .
$$
In polar coordinates, $\E^{b_\eta}(v)=\|F\|_{L^p(\mu)}^p$. Since $(1-s)\int_0^{\delta}r^{p(1-s)-1}dr=\frac1p\delta^{p(1-s)}$,
$$
\|G\|_{L^p(\mu)}^p=\delta^{p(1-s)}\intr\mathcal B_\eta\Big(\frac x\ve,\nabla v\Big)dx .
$$
Since $v(x)-v(x-r\omega)-r\nabla v(x)\cdot\omega=-r^2\int_0^1(1-t)D^2v(x-tr\omega)[\omega,\omega]\,dt$, we have $\|F(r,\omega,\cdot)-G(r,\omega,\cdot)\|_{L^p(dx)}\le\frac r2\|D^2v\|_p$, and by Lemma \ref{lema.der},
\begin{align*}
\|F-G\|_{L^p(\mu)}^p&\le m_+n\omega_n(1-s)\frac{\|D^2v\|_p^p}{2^p}\int_0^{\delta}r^{p(1-s)+p-1}dr\\
&\le C_1(1-s)\,\delta^{p(2-s)}\,\delta^{p(s-2)}\Lambda=C_1(1-s)\Lambda .
\end{align*}
The result follows from Minkowski's inequality in $L^p(\mu)$ and Steps 1 and 2.
\end{proof}

\begin{lema}\label{lema.conj}
Let $b_k\colon\R^n\times\Sn\to[m_-,m_+]$ be measurable with $b_k\to a$ a.e.\ in $\R^n\times\Sn$, and let $\mathcal B_k$ be the corresponding densities. Then, for a.e.\ $y\in\R^n$, $\mathcal B_k(y,\xi)\to\A(y,\xi)$ for every $\xi$ and $\mathcal B_k^c(y,\sigma)\to\A^c(y,\sigma)$ for every $\sigma$.
\end{lema}

\begin{proof}
By Fubini's theorem, for a.e.\ $y$ we have $b_k(y,\omega)\to a(y,\omega)$ for a.e.\ $\omega$, hence $\mathcal B_k(y,\xi)\to\A(y,\xi)$ for every $\xi$ by dominated convergence. Since these are convex functions, the convergence is locally uniform in $\xi$; by the uniform coercivity in \eqref{crec} the supremum defining the conjugates is attained in a ball depending only on $\sigma$, and therefore $\mathcal B_k^c(y,\sigma)\to\A^c(y,\sigma)$.
\end{proof}

\begin{lema}[Cell problem]\label{lema.celda}
For every $\xi_0\in\R^n$ there exists $\tau\in L^{p'}_{\rm loc}(\R^n;\R^n)$, $Q$-periodic, with $\operatorname{div}\tau=0$ in $\mathcal D'(\R^n)$, such that, writing $\bar\tau=\int_Q\tau$,
\begin{equation}\label{celda}
\int_Q\A^c(y,\tau(y))\,dy=\bar\tau\cdot\xi_0-\A^*(\xi_0),\qquad\A^*(\xi)\ge\A^*(\xi_0)+\bar\tau\cdot(\xi-\xi_0)\quad\text{for all }\xi\in\R^n,
\end{equation}
and $|\bar\tau|\le C|\xi_0|^{p-1}$ with $C$ depending only on $n,p,m_\pm$. Moreover, $\A^*$ is strictly convex.
\end{lema}

\begin{proof}
By \eqref{crec} and the direct method there is a minimizer $w\in W^{1,p}_{\rm per}(Q)$ of $v\mapsto\int_Q\A(y,\xi_0+\nabla v)$. Let $\tau(y)=\nabla_\xi\A(y,\xi_0+\nabla w(y))$, extended periodically. The Euler--Lagrange equation reads $\int_Q\tau\cdot\nabla\varphi=0$ for every $\varphi\in W^{1,p}_{\rm per}(Q)$. If $\psi\in C^\infty_c(\R^n)$, then $\Phi(y)=\sum_{j\in\mathbb Z^n}\psi(y+j)$ is smooth and periodic and $\intr\tau\cdot\nabla\psi=\int_Q\tau\cdot\nabla\Phi=0$, so $\operatorname{div}\tau=0$ in $\mathcal D'(\R^n)$. By the equality case in Fenchel's inequality, $\A^c(y,\tau)=\tau\cdot(\xi_0+\nabla w)-\A(y,\xi_0+\nabla w)$; integrating over $Q$ and using $\int_Q\tau\cdot\nabla w=0$ we get the first identity in \eqref{celda}. If $\xi\in\R^n$ and $w_\xi$ is a minimizer for $\xi$, Fenchel's inequality gives
$$
\A^*(\xi)=\int_Q\A(y,\xi+\nabla w_\xi)\ge\int_Q\tau\cdot(\xi+\nabla w_\xi)-\int_Q\A^c(y,\tau)=\bar\tau\cdot\xi-\bar\tau\cdot\xi_0+\A^*(\xi_0).
$$
By \eqref{crec}, $c_1\|\xi_0+\nabla w\|_{L^p(Q)}^p\le\A^*(\xi_0)\le c_2|\xi_0|^p$ and $\|\tau\|_{L^{p'}(Q)}\le m_+K_{n,p}\|\xi_0+\nabla w\|_{L^p(Q)}^{p-1}$, which gives the bound on $\bar\tau$. Finally, $\A(y,\cdot)$ is strictly convex (since $p>1$ and $a\ge m_->0$); if $\A^*(\frac{\xi+\eta}2)=\frac12(\A^*(\xi)+\A^*(\eta))$, testing with $\frac12(w_\xi+w_\eta)$ shows that equality holds a.e.\ in the convexity inequality for $\A(y,\cdot)$ at the points $\xi+\nabla w_\xi$ and $\eta+\nabla w_\eta$, so they coincide a.e., and integrating over $Q$ gives $\xi=\eta$.
\end{proof}

\begin{prop}[Homogenization lower bound with varying coefficients]\label{prop.homog}
Let $b_k\colon\R^n\times\Sn\to[m_-,m_+]$ be measurable, $Q$-periodic in $y$, with $b_k\to a$ a.e.\ in $Q\times\Sn$, and let $\mathcal B_k$ be the corresponding densities. Let $\ve_k\downarrow0$ and let $\{v_k\}$ be bounded in $W^{1,p}(\R^n)$ with $v_k\to u$ in $L^p_{\rm loc}(\R^n)$. Then $u\in W^{1,p}(\R^n)$ and
$$
\intr\A^*(\nabla u)\,dx\le\liminf_{k\to\infty}\intr\mathcal B_k\Big(\frac x{\ve_k},\nabla v_k\Big)dx .
$$
\end{prop}

\begin{proof}
Clearly $u\in W^{1,p}(\R^n)$. By Lemma \ref{lema.conj} (and periodicity), $\mathcal B_k^c(y,\sigma)\to\A^c(y,\sigma)$ for every $\sigma$, for a.e.\ $y$.

Fix $\xi_0\in\R^n$, let $\tau$ be given by Lemma \ref{lema.celda}, and let $g_k(y)=\mathcal B_k^c(y,\tau(y))$, $g(y)=\A^c(y,\tau(y))$. Then $g_k\to g$ a.e.\ and $0\le g_k\le c_3|\tau|^{p'}\in L^1(Q)$, so $g_k\to g$ in $L^1(Q)$. Let $\psi\in C^1_c(\R^n)$, $\psi\ge0$. By Fenchel's inequality,
$$
\intr\psi\,\mathcal B_k\Big(\frac x{\ve_k},\nabla v_k\Big)dx\ge\intr\psi\,\tau\Big(\frac x{\ve_k}\Big)\cdot\nabla v_k\,dx-\intr\psi\,g_k\Big(\frac x{\ve_k}\Big)dx .
$$
Since $\tau(\cdot/\ve_k)\in L^{p'}_{\rm loc}$ is divergence free and $\psi v_k\in W^{1,p}$ has compact support, the first integral equals $-\intr v_k\,\tau(\frac x{\ve_k})\cdot\nabla\psi\,dx$; since $v_k\to u$ in $L^p_{\rm loc}$ and $\tau(\cdot/\ve_k)\rightharpoonup\bar\tau$ weakly in $L^{p'}_{\rm loc}$, it converges to $-\intr u\,\bar\tau\cdot\nabla\psi=\intr\psi\,\bar\tau\cdot\nabla u$. For the second one, covering $\operatorname{supp}\psi$ by cubes of side $\ve_k$,
$$
\Big|\intr\psi\,(g_k-g)\Big(\frac x{\ve_k}\Big)dx\Big|\le C_\psi\|g_k-g\|_{L^1(Q)}\to0,\qquad\intr\psi\,g\Big(\frac x{\ve_k}\Big)dx\to\intr\psi\,dx\int_Qg\,dy .
$$
By \eqref{celda} we conclude that
\begin{equation}\label{afin}
\liminf_{k\to\infty}\intr\psi\,\mathcal B_k\Big(\frac x{\ve_k},\nabla v_k\Big)dx\ge\intr\psi\,\ell_{\xi_0}(\nabla u)\,dx,\qquad\ell_{\xi_0}(\zeta)=\A^*(\xi_0)+\bar\tau\cdot(\zeta-\xi_0).
\end{equation}
Let $L$ be the liminf in the statement; we may assume $L<\infty$ and pass to a subsequence along which the integrals converge to $L$ and the measures $\mu_k=\mathcal B_k(\frac x{\ve_k},\nabla v_k)\,dx$ converge weakly-$*$ to a Radon measure $\mu$, with $\mu(\R^n)\le L$. By \eqref{afin}, $\mu\ge\ell_{\xi_0}(\nabla u)\,dx$ for every $\xi_0$ in a countable dense set $D\subset\R^n$, so the density of the absolutely continuous part of $\mu$ is at least $\sup_{\xi_0\in D}\ell_{\xi_0}(\nabla u(x))$ a.e. By \eqref{celda}, $\ell_{\xi_0}\le\A^*$ with equality at $\xi_0$; since $\A^*$ is continuous and $|\bar\tau|\le C|\xi_0|^{p-1}$, $\sup_{\xi_0\in D}\ell_{\xi_0}(\zeta)=\A^*(\zeta)$ for every $\zeta$. Hence $L\ge\mu(\R^n)\ge\intr\A^*(\nabla u)$.
\end{proof}

\begin{proof}[Proof of Theorem \ref{teo.micro}]
Let $\Lambda_0=\sup_k\E_k(u_k)$, $\delta_k=\ve_k/R_k$ and $\eta_k=\delta_k/\ve_k=R_k^{-1}\to0$. Then $\delta_k\to0$ and $\delta_k^{p(1-s_k)}\to\theta>0$. Let $v_k=(u_k)_{\delta_k}$, and let $b_k=b_{\eta_k}$, $\mathcal B_k=\mathcal B_{\eta_k}$ be given by \eqref{beta}; they are $Q$-periodic in $y$ and $b_k\to a$ a.e. By Proposition \ref{prop.red} and \eqref{coerc}, for $k$ large,
$$
\E_k(u_k;0,\delta_k)^{1/p}\ge\Big(1-\frac{C_0\delta_k}{m_-}\Big)^{1/p}\Big[\Big(\delta_k^{p(1-s_k)}\intr\mathcal B_k\Big(\frac x{\ve_k},\nabla v_k\Big)dx\Big)^{1/p}-\big(C(1-s_k)\Lambda_0\big)^{1/p}\Big].
$$
By Lemma \ref{lema.moll}(1), $\|v_k-u_k\|_p\to0$, so $v_k\to u$ in $L^p_{\rm loc}$ and $\{v_k\}$ is bounded in $L^p$. By Lemma \ref{lema.der}, $\|\nabla v_k\|_p\le C\delta_k^{s_k-1}=C(\delta_k^{p(1-s_k)})^{-1/p}$, which is bounded because $\theta>0$. Hence $\{v_k\}$ is bounded in $W^{1,p}(\R^n)$, and Proposition \ref{prop.homog} gives $u\in W^{1,p}(\R^n)$ and $\liminf_k\intr\mathcal B_k(\frac x{\ve_k},\nabla v_k)\ge\intr\A^*(\nabla u)$. Since $\delta_k^{p(1-s_k)}\to\theta$, the result follows.
\end{proof}

\begin{rem}
The boundedness of $\nabla v_k$ in $L^p$ is exactly where $\theta>0$ is used: the mollification scale $\delta_k$ is far below $\ve_k$ in the linear scale, but at logarithmic distance $o(1/(1-s_k))$ from it.
\end{rem}

\subsection{The upper bound}\label{sec.limsup}

\begin{lema}[Large scales, upper bound]\label{lema.suave}
Assume \eqref{H1}, \eqref{H3} and \eqref{theta}. Then, for every $u\in W^{1,p}(\R^n)$,
$$
\limsup_{k\to\infty}\E_k(u;R_k\ve_k,\infty)\le(1-\theta)\intr\bar\A(\nabla u)\,dx .
$$
\end{lema}

\begin{proof}
Let $\rho_k=R_k\ve_k$. If $\theta=1$, Lemma \ref{lema.cota} gives
$$
\E_k(u;\rho_k,\infty)\le c_2\big(1-\min\{\rho_k,1\}^{p(1-s_k)}\big)\|\nabla u\|_p^p+C(1-s_k)\|u\|_p^p\to0 .
$$
Assume $\theta<1$, so that $\rho_k\to0$.

\emph{The case $\phi\in C^2_c(\R^n)$.} Fix $\kappa\in(0,1)$. First, $\E_k(\phi;\kappa,\infty)\le m_+(1-s_k)2^p\|\phi\|_p^p\frac{n\omega_n}{s_kp}\kappa^{-s_kp}\to0$. Next, by the second inequality in Lemma \ref{lema.prom} with $v=\phi$, and since $\|\tau_{-h}\phi-\phi\|_p\le|h|\|\nabla\phi\|_p$ and $\|\tau_{-h}\nabla\phi-\nabla\phi\|_p\le2\|\nabla\phi\|_p$,
$$
\big|\E_k(\phi;\rho_k,\kappa)-\bar\E_{s_k}(\phi;\rho_k,\kappa)\big|\le C_\phi\,\ve_k\,(1-s_k)\int_{\rho_k}^\kappa r^{p(1-s_k)-2}\,dr\le C_\phi\,(1-s_k)\,\frac{\ve_k}{\rho_k}\to0 ,
$$
where we used $r^{p(1-s_k)}\le1$ and $\int_{\rho_k}^\kappa r^{-2}dr\le\rho_k^{-1}$. Finally, arguing as in the proof of Lemma \ref{lema.banda} (now with the upper bounds $\alpha^p\le\beta^p+Cr$ and $\bar m(r\omega)\le\bar a(\omega)+C_0r$, and integrating over a bounded set containing the $1$-neighborhood of $\operatorname{supp}\phi$),
$$
\bar\E_{s_k}(\phi;\rho_k,\kappa)\le\big(\kappa^{p(1-s_k)}-\rho_k^{p(1-s_k)}\big)\intr\bar\A(\nabla\phi)\,dx+C_\phi(1-s_k)\to(1-\theta)\intr\bar\A(\nabla\phi)\,dx .
$$

\emph{The general case.} By Lemma \ref{lema.cota}, $\E_k(\cdot\,;\rho_k,\infty)^{1/p}$ is a seminorm bounded by $C\|\cdot\|_{W^{1,p}}$ uniformly in $k$. Hence, for $\phi\in C^2_c$,
$$
\limsup_k\E_k(u;\rho_k,\infty)^{1/p}\le\Big((1-\theta)\intr\bar\A(\nabla\phi)\Big)^{1/p}+C\|u-\phi\|_{W^{1,p}},
$$
and the result follows by density, since $v\mapsto\intr\bar\A(\nabla v)$ is continuous in $W^{1,p}$.
\end{proof}

\begin{prop}[Recovery sequence]\label{prop.rec}
Assume \eqref{H1}, \eqref{H3} and \eqref{theta}. Let $u$ be continuous, piecewise affine and with compact support. Then for every $\gamma>0$ there exist Lipschitz functions $u_k$, with $\operatorname{supp}u_k\subset\operatorname{supp}u$ and $\|u_k-u\|_\infty\le C\ve_k$, such that
$$
\limsup_{k\to\infty}\E_k(u_k)\le\E_\theta(u)+\gamma .
$$
\end{prop}

\begin{proof}
Let $\{T_i\}_{i=1}^N$ be closed simplices with disjoint interiors such that $\operatorname{supp}u=\bigcup_iT_i=:K$ and $\nabla u=\xi_i$ in the interior of $T_i$. The parameters are chosen in the following order. First, for each $i$ choose $w_i\in C^\infty_{\rm per}(Q)$ such that
$$
\int_Q\A(y,\xi_i+\nabla w_i(y))\,dy\le\A^*(\xi_i)+\gamma,
$$
which is possible by the density of smooth periodic functions in $W^{1,p}_{\rm per}(Q)$ and \eqref{crec}, and let $M_0=1+\max_i(|\xi_i|+\|\nabla w_i\|_\infty)$. Next, let $\gamma'>0$ be such that $m_+M_0^pn\omega_n\gamma'|K|\le\gamma$, let $T_i'\subset\operatorname{int}T_i$ be a closed simplex homothetic to $T_i$ with $|T_i\setminus T_i'|\le\gamma'|T_i|$, and let $\varphi_i\in C^\infty_c(\operatorname{int}T_i)$ with $0\le\varphi_i\le1$ and $\varphi_i=1$ on $T_i'$. Define
$$
u_k(x)=u(x)+\ve_k\sum_{i=1}^N\varphi_i(x)\,w_i\Big(\frac x{\ve_k}\Big).
$$
Then $\operatorname{supp}u_k\subset K$, $\|u_k-u\|_\infty\le C\ve_k$ and, since $\nabla u_k=\xi_i+\varphi_i\nabla w_i(\frac x{\ve_k})+\ve_kw_i(\frac x{\ve_k})\nabla\varphi_i$ in $T_i$, we have $\|\nabla u_k\|_\infty\le M_0$ for $k$ large. Let $\rho_k^-=\ve_k/R_k$, $\rho_k^+=R_k\ve_k$, and split $\E_k(u_k)$ according to the bands $|h|<\rho_k^-$, $\rho_k^-\le|h|\le\rho_k^+$ and $|h|>\rho_k^+$. In all cases the integrand vanishes unless $x\in K$ or $x-h\in K$.

\emph{Intermediate scales.} Since $|u_k(x)-u_k(x-h)|\le M_0|h|$, and for each $h$ the set of $x$ where the integrand does not vanish has measure at most $2|K|$,
$$
\E_k(u_k;\rho_k^-,\rho_k^+)\le2m_+M_0^p|K|\,n\omega_n(1-s_k)\int_{\rho_k^-}^{\rho_k^+}r^{p(1-s_k)-1}dr=C\big((\rho_k^+)^{p(1-s_k)}-(\rho_k^-)^{p(1-s_k)}\big)\to0 ,
$$
because both terms converge to $\theta$.

\emph{Large scales.} Let $g_k=u_k-u$, so that $\|g_k\|_p\le C\ve_k$. By Lemma \ref{lema.cota},
$$
\E_k(g_k;\rho_k^+,\infty)\le C(1-s_k)(\rho_k^+)^{-s_kp}\|g_k\|_p^p\le C(1-s_k)\,\ve_k^{p(1-s_k)}R_k^{-s_kp}\to0 .
$$
Since $\E_k(\cdot\,;\rho_k^+,\infty)^{1/p}$ is a seminorm, Lemma \ref{lema.suave} gives
$$
\limsup_k\E_k(u_k;\rho_k^+,\infty)=\limsup_k\E_k(u;\rho_k^+,\infty)\le(1-\theta)\intr\bar\A(\nabla u)\,dx .
$$

\emph{Small scales.} Let $G_k^i=\{x\in T_i'\colon\operatorname{dist}(x,\partial T_i')>\rho_k^-\}$ and $N_k=(K+B_{\rho_k^-})\setminus\bigcup_iG_k^i$. Then $|N_k|\to\sum_i|T_i\setminus T_i'|\le\gamma'|K|$. If $|h|<\rho_k^-$ and $x\notin K+B_{\rho_k^-}$ the integrand vanishes. On $N_k$ we use the Lipschitz bound:
$$
(1-s_k)\int_{\{|h|<\rho_k^-\}}\int_{N_k}m\Big(\frac x{\ve_k},h\Big)\frac{|u_k(x)-u_k(x-h)|^p}{|h|^{n+s_kp}}\,dx\,dh\le\frac{m_+M_0^pn\omega_n}p|N_k|,
$$
whose limsup is at most $\gamma$ by the choice of $\gamma'$. If $x\in G_k^i$ and $|h|<\rho_k^-$, then the segment $[x-h,x]$ is contained in $T_i'$, where $u_k(z)=u(z)+\ve_kw_i(z/\ve_k)$ and $u$ is affine with gradient $\xi_i$. Hence, with $h=r\omega$,
$$
\Big|\frac{u_k(x)-u_k(x-h)}r-\Big(\xi_i+\nabla w_i\Big(\frac x{\ve_k}\Big)\Big)\cdot\omega\Big|\le\frac{r}{2\ve_k}\|D^2w_i\|_\infty\le\frac{C}{R_k},
$$
and both quantities are bounded by $M_0$, so their $p$-th powers differ by at most $C/R_k$. Using also $m(y,r\omega)\le a(y,\omega)+C_0r$ and $(1-s_k)\int_0^{\rho_k^-}r^{p(1-s_k)-1}dr=\frac1p(\rho_k^-)^{p(1-s_k)}\le\frac1p$,
\begin{align*}
(1-s_k)\int_{\{|h|<\rho_k^-\}}\int_{G_k^i}m\Big(\frac x{\ve_k},h\Big)&\frac{|u_k(x)-u_k(x-h)|^p}{|h|^{n+s_kp}}\,dx\,dh\\
&\le(\rho_k^-)^{p(1-s_k)}\int_{T_i}\A\Big(\frac x{\ve_k},\xi_i+\nabla w_i\Big(\frac x{\ve_k}\Big)\Big)dx+C\Big(\frac1{R_k}+\rho_k^-\Big).
\end{align*}
Since $y\mapsto\A(y,\xi_i+\nabla w_i(y))$ is bounded and periodic, the last integral converges to
$$
|T_i|\int_Q\A(y,\xi_i+\nabla w_i)\,dy\le|T_i|\big(\A^*(\xi_i)+\gamma\big).
$$
Adding up,
$$
\limsup_k\E_k(u_k;0,\rho_k^-)\le\theta\intr\A^*(\nabla u)\,dx+\gamma\,(|K|+1) .
$$

Adding the three bands we get $\limsup_k\E_k(u_k)\le\E_\theta(u)+\gamma(|K|+1)$; applying this with $\gamma/(|K|+1)$ in place of $\gamma$ gives the result.
\end{proof}

\subsection{Proofs of Theorem \ref{teo.main} and Corollary \ref{cor.min}}\label{sec.pruebas}

\begin{proof}[Proof of Theorem \ref{teo.main}]
(1) Let $\Lambda_0=\sup_k\E_k(u_k)$. By \eqref{coerc} and Lemma \ref{lema.transl}, $\|\tau_hu_k-u_k\|_p^p\le C_{n,p}|h|^{s_kp}\Lambda_0/m_-$ for every $h$ and $k$. Since $s_k\ge\frac12$ for $k$ large, the sequence is bounded and equicontinuous in $L^p$ with respect to translations, and the Riesz--Fr\'echet--Kolmogorov theorem gives the precompactness in $L^p_{\rm loc}(\R^n)$. If $u$ is a limit point, Fatou's lemma gives $\|\tau_hu-u\|_p^p\le C_{n,p}|h|^p\Lambda_0/m_-$ for every $h$, and since $p>1$ this implies $u\in W^{1,p}(\R^n)$.

(2) We may assume that $\liminf_k\E_k(u_k)<\infty$ and, passing to a subsequence (which still satisfies \eqref{theta}), that $\sup_k\E_k(u_k)<\infty$. Since the bands $\{|h|<\ve_k/R_k\}$ and $\{|h|>R_k\ve_k\}$ are disjoint and the integrand is nonnegative,
$$
\E_k(u_k)\ge\E_k(u_k;0,\ve_k/R_k)+\E_k(u_k;R_k\ve_k,\infty).
$$
If $0<\theta<1$, Theorems \ref{teo.micro} and \ref{teo.macro} give $u\in W^{1,p}$ and $\liminf_k\E_k(u_k)\ge\theta\intr\A^*(\nabla u)+(1-\theta)\intr\bar\A(\nabla u)$. If $\theta=0$ we use only Theorem \ref{teo.macro}, and if $\theta=1$ only Theorem \ref{teo.micro}.

(3) Let $\E''(u)=\inf\{\limsup_k\E_k(u_k)\colon u_k\to u\text{ in }L^p(\R^n)\}$ be the $\Gamma$-upper limit, which is lower semicontinuous in $L^p(\R^n)$ \cite[Proposition 6.8]{DalMaso}; moreover, the infimum in its definition is attained \cite[Proposition 8.1]{DalMaso}, so that, by a diagonal argument, recovery sequences can be taken among the functions constructed in Proposition \ref{prop.rec} for a sequence of approximants. By Proposition \ref{prop.rec}, $\E''(u)\le\E_\theta(u)$ for every continuous, piecewise affine $u$ with compact support. These functions are dense in $W^{1,p}(\R^n)$ and $\E_\theta$ is continuous in $W^{1,p}(\R^n)$ by \eqref{growth}, so $\E''\le\E_\theta$ in $W^{1,p}(\R^n)$. If $\Omega$ is bounded and Lipschitz and $u=0$ a.e.\ outside $\Omega$, then $u|_\Omega\in W^{1,p}_0(\Omega)$ and $u$ is the limit in $W^{1,p}$ of piecewise affine functions with compact support in $\Omega$ (approximate first by functions in $C^\infty_c(\Omega)$ and then interpolate on a fine triangulation); the same argument, with $\E''$ defined through sequences vanishing outside $\Omega$, applies, and the recovery sequences of Proposition \ref{prop.rec} are supported in $\Omega$.
\end{proof}

\begin{proof}[Proof of Corollary \ref{cor.min}]
Let $d=\operatorname{diam}\Omega+1$. If $v=0$ a.e.\ outside $\Omega$ and $|h|=d$, then $\tau_hv$ and $v$ have disjoint supports, and Lemma \ref{lema.transl} and \eqref{coerc} give the Poincar\'e inequality $2\|v\|_p^p\le C_{n,p}d^{sp}\E_{\ve,s}(v)/m_-$, uniform in $\ve$ and $s$. Hence the functionals $\E_k(v)-\int_\Omega fv$ are coercive and strictly convex on $\{v=0\text{ a.e.\ outside }\Omega\}$, and lower semicontinuous in $L^p$ by Fatou's lemma; since they are convex, they are also weakly lower semicontinuous in the reflexive space $L^p(\Omega)$, and this gives existence and uniqueness of $u_k$. Comparing with $v=0$, $\E_k(u_k)\le\|f\|_{p'}\|u_k\|_p\le C\E_k(u_k)^{1/p}$, so $\sup_k\big(\E_k(u_k)+\|u_k\|_p\big)<\infty$ because $p>1$. By Theorem \ref{teo.main}(1), and since the $u_k$ vanish outside $\Omega$, $\{u_k\}$ is precompact in $L^p(\Omega)$. Since $v\mapsto\int_\Omega fv$ is continuous in $L^p$, Theorem \ref{teo.main} and the fundamental theorem of $\Gamma$-convergence give that every limit point minimizes the limit problem. Its minimizer is unique because $\A_\theta$ is strictly convex (Lemma \ref{lema.celda} for $\A^*$; for $\bar\A$ it is immediate).
\end{proof}

\section{Further comments}\label{sec.com}

\begin{itemize}
\item \emph{Symmetry.} The kernel $m(\frac x\ve,h)$ does not inherit the symmetry $m(x,h)=m(x-h,-h)$ assumed in \cite{FBS-Israel} to write the Euler--Lagrange operator; the operator associated to $\E_{\ve,s}$ must be written with the symmetrized kernel.
\item \emph{Other kernels.} For kernels oscillating also in $h$ at scale $\ve$, such as $m(\frac x\ve,\frac h\ve)$, $c(\frac x\ve)c(\frac{x-h}\ve)$ or the coefficients $k(\frac x\ve,\frac y\ve)$ of \cite{BMS}, we expect that the large and the small scales can be treated in a similar way, with different limits of the kernels in each range: for $k(\frac x\ve,\frac y\ve)$ with $k$ continuous, the small scales see the diagonal values $k(y,y)$ and the large scales the average of $a$ in both variables, as proved in \cite{BMS} for $p=2$.
\item \emph{Non periodic coefficients.} The proof for the large scales only uses averaging, and the small scales reduce to local homogenization. Stochastic or general $G$-convergent structures should work similarly.
\end{itemize}

\end{document}